\documentclass[12pt]{amsart}
 \usepackage[latin1]{inputenc}
 \usepackage[dvips]{graphicx}
 \usepackage{wrapfig}
 \usepackage{amsmath}
 \usepackage{amsthm}
 \usepackage{amsfonts}
 \usepackage{amssymb}
 \usepackage{layout}
 \usepackage{verbatim}
 \usepackage{alltt}
\usepackage{yfonts}
\usepackage{stmaryrd}
 \usepackage{tipa}
 \usepackage{indentfirst}

\usepackage[all]{xy}
\usepackage{setspace}

\newtheorem*{thma}{Theorem A}

\newtheorem*{conj}{Conjecture}

\newcounter{thm}[section]

\renewcommand{\thethm} {\arabic{section}.\arabic{thm}}

\newenvironment{thm}{\par\medskip\noindent\refstepcounter{thm}
\bgroup{\hspace*{-0.15 cm}\bf{Theorem}
\thethm.}\bgroup\it}{\egroup \egroup\par\medskip}

\newenvironment{lemma}{\par\medskip\noindent\refstepcounter{thm}
\bgroup{\hspace*{-0.15 cm}\bf{Lemma} \thethm.}\bgroup\it}{\egroup
\egroup\par\medskip}

\newenvironment{prop}{\par\medskip\noindent\refstepcounter{thm}
\bgroup{\hspace*{-0.15 cm}\bf{Proposition}
\thethm.}\bgroup\it}{\egroup \egroup\par\medskip}

\newenvironment{cor}{\par\medskip\noindent\refstepcounter{thm}
\bgroup{\hspace*{-0.15 cm}\bf{Corollary}
\thethm.}\bgroup\it}{\egroup \egroup\par\medskip}

\newenvironment{define}{\par\medskip\noindent\refstepcounter{thm}
\bgroup{\hspace*{-0.15 cm}\bf{Definition}
\thethm.}\bgroup}{\egroup \egroup\par\medskip}

\newenvironment{example}{\par\medskip\noindent\refstepcounter{thm}
\bgroup{\hspace*{-0.15 cm}\bf{Example}
\thethm.}\bgroup}{\egroup \egroup\par\medskip}

\newenvironment{rem}{\par\medskip\noindent\refstepcounter{thm}
\bgroup{\hspace*{-0.15 cm}\bf{Remark} \thethm.}\bgroup}{\egroup
\egroup\par\medskip} 

\newenvironment{conjecture}{\par\medskip\noindent\refstepcounter{thm}
\bgroup{\hspace*{-0.15 cm}\bf{Conjecture}
\thethm.}\bgroup\it}{\egroup \egroup\par\medskip}
\newcounter{Athm}[section]

\renewcommand{\theAthm} {\arabic{Athm}}

\newenvironment{Athm}{\par\medskip\noindent\refstepcounter{Athm}
\bgroup{\hspace*{-0.15 cm}\bf{Theorem}
A.\theAthm.}\bgroup\it}{\egroup \egroup\par\medskip}

\newenvironment{Aprop}{\par\medskip\noindent\refstepcounter{Athm}
\bgroup{\hspace*{-0.15 cm}\bf{Proposition}
A.\theAthm.}\bgroup\it}{\egroup \egroup\par\medskip}

\newenvironment{Adefine}{\par\medskip\noindent\refstepcounter{Athm}
\bgroup{\hspace*{-0.15 cm}\bf{Definition}
A.\theAthm.}\bgroup}{\egroup \egroup\par\medskip}

\newenvironment{Arem}{\par\medskip\noindent\refstepcounter{Athm}
\bgroup{\hspace*{-0.15 cm}\bf{Remark} A.\theAthm.}\bgroup}{\egroup
\egroup\par\medskip}

\newcommand{\LL}{\Lambda}
\newcommand{\QQ}{\mathbb{Q}}
\newcommand{\FF}{\mathcal{F}}
\newcommand{\FFc}{\mathcal{F}_{\textup{\lowercase{can}}}}
\newcommand{\lra}{\longrightarrow}
\newcommand{\ZZ}{\mathbb{Z}}
\newcommand{\PP}{\mathcal{P}}
\newcommand{\NN}{\mathcal{N}}
\newcommand{\ra}{\rightarrow}	
\newcommand{\xx}{\mathbf{X}}
\newcommand{\be}{\begin{equation}}
\newcommand{\ee}{\end{equation}}

\newcommand{\XX}{\mathcal{X}}

\newcommand{\qq}{\hbox{\frakfamily q}}
\newcommand{\KS}{\textbf{\textup{KS}}}

\newcommand{\tk}{\tilde{K}}
\newcommand{\tH}{\tilde{H}}
\newcommand{\tgamma}{\tilde{\gamma}}
\newcommand{\tlambda}{\tilde{\Lambda}}
\newcommand{\tGamma}{\tilde{\Gamma}}
\newcommand{\fontaineO}{\mathcal{O}_{\widehat{\varepsilon^{\textup{nr}}}}}
\newcommand{\fontaineOK}{\mathcal{O}_{\varepsilon(K)}}
\newcommand{\htam}{$\mathbf{\mathbb{H}.T}$}
\newcommand{\hsez}{$\mathbf{\mathbb{H}.sEZ}$}

\newcommand{\Gal}{\textup{Gal}}

\newcommand{\hone}{$\mathbf{H.1}$}

\newcommand{\hfour}{$\mathbf{H.4}$}

\newcommand{\LLL}{\mathbb{L}_{\infty}}
\newcommand{\hsezF}{$\mathbf{\mathbb{H}.sEZ}_{/F}$}
\newcommand{\htamF}{$\mathbf{\mathbb{H}.T}_{/F}$}
\newcommand{\Ll}{\mathcal{L}}
\newcommand{\kk}{\mathcal{K}}
\newcommand{\eee}{\epsilon}
\newcommand{\vv}{\mathbb{V}}
\newcommand{\al}{\mathcal{L}}
\newcommand{\hli}{\hbox{H}_{\all}}
\newcommand{\hhh}{\bigwedge^{r-1}\textbf{\textup{Hom}}(\mathbb{V},\LL\otimes_{\ZZ_p}\ZZ_p[[\textup{Gal}(M/k)]])}

\newcommand{\all}{\mathbb{L}}

\begin{document}
\title{Stark units and main conjectures for totally real fields}

\author{K\^az\i m B\"uy\"ukboduk}

\address{Kazim Buyukboduk \hfill\break\indent
Department of Mathematics
\hfill\break\indent 450 Serra Mall, Bldg 380 \hfill\break\indent Stanford, CA, 94305
\hfill\break\indent USA} \email{{\tt kazim@math.stanford.edu }\hfill\break\indent {\it
Web page:} {\tt math.stanford.edu/$\sim$kazim}}


\begin{abstract}

Main theorem of~\cite{kbb}  suggests that it should be possible to lift the Kolyvagin systems of Stark units constructed in~\cite{kbbstark} to a Kolyvagin system over the cyclotomic Iwasawa algebra. This is what we prove in this paper. This construction gives the first example towards a more systematic study of Kolyvagin system theory over an Iwasawa algebra when the core Selmer rank is greater than one.  As a result of this construction, we reduce the main conjectures of Iwasawa theory for totally real fields to a statement of local Iwasawa theory. This statement, however, turns out to be interesting in its own right as it suggests a relation between solutions to $p$-adic  and complex Stark conjectures.

\end{abstract}

\maketitle
\tableofcontents

\section{Introduction}
This is an attempt to understand the Kolyvagin system theory (over an Iwasawa algebra) when the core Selmer rank (in the sense of~\cite{mr02} Definitions 4.1.8 and 4.1.11) is greater than one, and should be considered as a continuation of our earlier paper~\cite{kbbstark}.

Kolyvagin system machinery is designed to bound the size of a Selmer group. In all well-known cases this bound obtained relates to $L$-values, and thus provides a link between arithmetic data to analytic data. Well-known prototypes for such a relation between arithmetic and analytic data are Birch and Swinnerton-Dyer conjecture (more generally Bloch-Kato conjectures) and main conjectures of Iwasawa theory. Kolyvagin system machinery has been successfully applied by  many to obtain deep results towards proving these conjectures.

In this paper, we concentrate on Kolyvagin systems over the cyclotomic Iwasawa algebra (which we henceforth denote by $\LL$). When the core Selmer rank is one, Kolyvagin systems over $\LL$ are proved to exist in wide variety of cases in~\cite{kbb}. Further, Mazur and Rubin shows in \cite{mr02} \S5.3 that these Kolyvagin systems can be used to compute the correct size of an appropriate Selmer group. However, when the core Selmer rank is greater than one not much is known.

The most basic  example of a core Selmer rank greater than one situation arises if one attempts to utilize the Euler system that would come from the Stark elements (whose existence were predicted by Rubin~\cite{ru96}).  Rubin was first to study the Euler system of Stark units in~\cite{ru92} where  he proved a Gras-type formula for the $\chi$-isotypic component of a certain ideal class group under certain assumptions on the character $\chi$ (which essentially ensured that the Selmer core rank of $T=\ZZ_p(1)\otimes\chi^{-1}$, in the sense of Definitions 4.8 and 4.1.11 of~\cite{mr02}, is one). These assumptions have been removed and a general Gras-type conjecture was proved in~\cite{kbbstark}. The proof relied on constructing auxiliary  Selmer structures in a systematic way to cut down the core Selmer rank to one, and using which one obtains a useful collection of Kolyvagin systems.

In this paper we basically lift the Kolyvagin systems for the modified Selmer structures  obtained  in~\cite{kbbstark} to Kolyvagin systems over the cyclotomic Iwasawa algebra. To be able to do this, we first modify (in \S \ref{sec:modified selmer}) the \emph{classical} Selmer structure along the cyclotomic tower. The main theorem of~\cite{kbb} shows that there are Kolyvagin systems over $\LL$ for the modified Selmer structure (see also \S\ref{kolsys-1} below). In \S\ref{kolsys-starkunits} we show how to obtain these Kolyvagin systems using the \emph{Euler system of Stark units} (which were introduced by Rubin~\cite{ru96}). In fact, the classical way of obtaining Kolyvagin systems is to apply Kolyvagin's decent on an Euler system. This has been systematized by Mazur and Rubin, who constructs what they call the \emph{Euler system to Kolyvagin system map} (\emph{c.f.} \cite{mr02} Theorem 5.3.3). The problem now is that if one applies the Euler system to Kolyvagin system map directly on the Euler system of Stark units, one does not in general get a Kolyvagin system for the modified Selmer structure. This is what we overcome in \S\ref{sec:homs}. Once we have the Kolyvagin system for the  modified Selmer structures, we apply the Kolyvagin system machinery of~\cite{mr02}.

Before we state our main results, we introduce some notation and set the hypotheses that will be assumed throughout this paper. 
\subsection{Notation}
\label{notation}
Throughout this paper the following notation and hypotheses are in effect.

$k$ is a totally real number field  of degree $r$ over $\QQ$. Fix once and for all an algebraic closure $\overline{k}$ of $k$, and a rational odd prime $p$. $k_\infty$ will denote the cyclotomic $\ZZ_p$-extension of $k$, and $G_k:=\Gal(\overline{k}/k)$, the absolute Galois group of $k$.

Let $\chi$ be a $totally\, even$ character of Gal$(\overline{k}/k)$ ($i.e.$ it is trivial on all complex conjugations inside Gal$(\overline{k}/k)$) into $\ZZ_p^{\times}$ that has finite order, and let $L$ be the fixed field of $\hbox{ker}(\chi)$ inside $\overline{k}$.  $f_{\chi}$ is the conductor of $\chi$, and $\Delta$ is the Galois group $\hbox{Gal}(L/k)$ of the extension $L/k$.

For any abelian group $A$, $A^{\wedge}$ will always denote its $p$-adic completion. Define also $A^{\chi}$ to be the $\chi$-isotypic component of $A^{\wedge}$.

$T$ will always stand for the  Gal$(\overline{k}/k)$-representation $\ZZ_p(1)\otimes~\chi^{-1}$ (except in the Appendix, where $T$ will be an arbitrary unramified almost everywhere $\ZZ_p [[\textup{Gal}(\overline{k}/k)]]$-module which is free of finite rank over $\ZZ_p$). 

For some results we will assume Leopoldt's conjecture. The following hypotheses will be occasionally assumed as well:
\begin{enumerate}
\item $p\nmid f_{\chi}$ (i.e. $L/k$ is unramified at all primes of $k$ above $p$),
\item $k$ is unramified at all primes above $p$.
\item $\chi(\textup{Frob}_\wp)\neq1$ for any prime $\wp$ of $k$ above $p$, where $\textup{Frob}_\wp$ denotes a Frobenius  element at $\wp$ inside $\Gal(\overline{k}/k)$. 
\end{enumerate}

(1) was already assumed in~\cite{kbbstark}, but we think that it could possibly be removed (both in~\cite{kbbstark} and here in this paper) using an argument similar to what we utilize to prove Proposition~\ref{big-free}. 

We are almost certain that this paper could have been written without (2) as well, but that would have at least necessitated fixing a prime $\wp$ of $k$ above $p$ to define our auxiliary local conditions (\S\ref{modify-at-p}), which we decided not to for the sake of consistency with~\cite{kbbstark}. In fact one of the major points we would like to study here is the $\LL$-adic Kolyvagin systems (in the sense of~\cite{kbb}) when core Selmer rank is strictly larger than \emph{one}. From this view point, we decided that we would rather assume (2) than losing our liberty in modifying the local conditions at $p$.

(3), however, is a more serious assumption. This is the assumption $\mathbf{\mathbb{H}.sEZ}$ of~\cite{kbb} \S2.2 translated to our setting; and it appears here for exactly the same reason why it appears in~\cite{kbb}. Note that, by assuming (3), we are actually requiring that  the Deligne-Ribet  $p$-adic $L$-function attached to $\chi$ has no $exceptional \, zeros$, in the sense of~\cite{g2, mtt}. 

\subsection{Statement of main results}
Suppose the set of places $S$ of $k$ contains no non-archimedean primes which split completely in $L/k$. Let $M_\infty$ is the maximal abelian $p$-extension of $L_\infty$ unramified outside primes above $p$, and set $\mathbb{T}=T\otimes\LL$. Let $H^1(k_p,\mathbb{T})$ denote the semi-local cohomology group\footnote{Let $k_n$ be the unique sub-extension of $k_\infty/k$ such that $[k_n:k]=p^n$. Set $L_n=L\cdot k_n$, and $\mathcal{U}_n^\chi$ denote the local units inside $(L_n\otimes\QQ_p)^\chi$. By Kummer theory (\emph{c.f.} \cite{r00} \S1.6.C, Proposition 3.2.6 and Lemma~\ref{at p: relaxed=unr} below) $H^1(k_p,\mathbb{T})$ may be identified with $\varprojlim_n \mathcal{U}_n^\chi$.} at $p$. Let $\mathbf{c}_{k_{\infty}}^{\textup{stark}}:=\{\varepsilon_{k_n}^\chi\}_n$, where $\varepsilon_{k_n}^\chi$ denotes an appropriate twist of Stark element (see \S\ref{sec:stark-euler} and \S\ref{sec:twisting} for details). For a torsion $\LL$-module $\mathbb{A}$, let $\textup{char}(\mathbb{A})$ denote its characteristic ideal.

\begin{thma}
$$\textup{char}\left(\Gal(M_\infty/L_\infty)^{\chi}\right)=\textup{char}\left(\bigwedge^rH^1(k_p,\mathbb{T})/\LL\cdot\mathbf{c}_{k_{\infty}}^{\textup{stark}}\right).$$
\end{thma}

When $k=\QQ$ (i.e. when $r=1$) $\mathbf{c}_{k_{\infty}}^{\textup{stark}}$ is given by the cyclotomic units. Further, the ideal on the right hand side of above equality in the statement of Theorem A is generated by the Kubota-Leopoldt $p$-adic $L$-function. This fact goes back to Iwasawa \cite{iwasawa-1}. Therefore, when $k=\QQ$, Theorem A is equivalent to the main conjectures of Iwasawa theory. 

Let $\mathcal{L}_{k}^{\chi}$ denote the Deligne-Ribet $p$-adic $L$-function attached to the character $\chi$ (see~\cite{deligne-ribet} for the construction of this $p$-adic $L$-function). Motivated by above example (when $k=\QQ$) and in the spirit of~ \cite{pr-kubota}, \cite{pr2} we propose the following:
\begin{conj}
$\mathcal{L}_{k}^{\chi}$ generates $\textup{char}\left(\bigwedge^rH^1(k_p,\mathbb{T})/\LL\cdot\mathbf{c}_{k_{\infty}}^{\textup{stark}}\right)$.
\end{conj}

It seems more reasonable to expect that such a relation, as we propose in above conjecture, should rather exist between the Deligne-Ribet $p$-adic $L$-function and ''\emph{p}-adic'' Stark elements (which would rather be solutions to a $p$-adic Stark conjecture). On the other hand, this conjecture is equivalent to the main conjectures of Iwasawa theory for totally real number fields, in particular above conjecture is true by~\cite{wiles-mainconj}. This whole discussion suggests a link between solutions to $p$-adic Stark conjectures and complex Stark conjectures at $s=0$. In a future paper, we hope to formulate this relation more precisely and hopefully prove the conjecture above.

\subsection{Acknowledgements} I would like to thank David Solomon, Christian Popescu and Karl Rubin for their comments and suggestions on this project.
\newpage
\section{Modified Selmer structures}
\label{sec:modified selmer}
\subsection{Selmer groups for $T=\ZZ_p(1)\otimes\chi^{-1}$}
Below we use the notation that was set in \S\ref{notation}. Let $\Gamma:=\Gal(k_\infty/k)$, and $\LL:=\ZZ_p[[\Gamma]]$ be the cyclotomic Iwasawa algebra. 

We first recall Mazur and Rubin's definition of a \emph{Selmer structure}, in particular the \emph{canonical Selmer structure} on $T\otimes\LL$.

\subsubsection{Local conditions}
\label{local conditions}
 Let $R$ be a complete local noetherian ring, and let $M$ be a $R[[G_k]]$-module which is free of finite rank over $R$. In this paper we will only be interested in the case when $R=\LL$ or its certain quotients, and $M$ is $T\otimes\LL$ or its relevant quotients by an ideal of $\LL$.     

For each place $\lambda$ of $k$, a \emph{local condition} $\FF$ (at $\lambda$)  on $M$ is a choice of an $R$-submodule $H^1_{\FF}(k_\lambda,M)$ of $H^1(k_\lambda,M)$. For the prime $p$, a local condition $\FF$ at $p$ will be a choice of an $R$-submodule $H^1_{\FF}(k_p,M)$ of the semi-local cohomology group $H^1(k_p,M):=\oplus_{\wp|p}H^1(k_\wp,M)$, where the direct sum is over all the primes $\wp$ of $k$ which lie above $p$.

For examples of local conditions see~\cite{mr02} Definitions 1.1.6 and 3.2.1.

Suppose that $\FF$ is a local condition (at the prime $\lambda$ of $k$)  on $M$. If $M^{\prime}$ is a submodule of $M$ (\emph{resp.} $M^{\prime \prime}$ is a quotient module), then $\FF$ induces local conditions (which we still denote by $\FF$) on $M^{\prime}$ (\emph{resp.} on $M^{\prime \prime}$), by taking $H^1_{\FF}(k_\lambda,M^{\prime})$ (\emph{resp.} $H^1_{\FF}(k_\lambda,M^{\prime \prime})$) to be the inverse image (\emph{resp.} the image) of $H^1_{\FF}(k_\lambda,M)$ under the natural maps induced by $$M^{\prime} \hookrightarrow M, \,\,\, \,\,\,\,\,\,\, M \twoheadrightarrow M^{\prime \prime}.$$

\begin{define}
\label{def:propagation}
\emph{Propagation} of a local condition $\FF$ on $M$ to a submodule $M^{\prime}$ (and a quotient $M^{\prime \prime}$ of $M$ is the local condition $\FF$ on $M^{\prime}$ (and on $M^{\prime \prime}$) obtained following the above procedure.
\end{define}
For example, if $I$ is an ideal of $R$, then a local condition on $M$ induces local conditions on $M/IM$ and $M[I]$, by \emph{propagation}.

\begin{define}
\label{cartier dual}
Define the \emph{Cartier dual} of $M$ to be the $R[[G_k]]$-module $$M^*:=\textup{Hom}(M,\mu_{p^{\infty}})$$ where $\mu_{p^{\infty}}$ stands for the $p$-power roots of unity.
\end{define}
Let $\lambda$ be a prime of $k$. There is the perfect local Tate pairing $$<\,,\,>_\lambda\,:H^1(k_\lambda,M) \times H^1(k_\lambda,M^*) \lra H^2(k_\lambda,\mu_{p^{\infty}}) \stackrel{\sim}{\lra}\QQ_p/\ZZ_p$$
\begin{define}
\label{dual local condition}
The \emph{dual local condition} $\FF^*$ on $M^*$ of a local  condition $\FF$ on $M$ is defined so that $H^1_{\FF^*}(k_\lambda,M^*)$ is the orthogonal complement of $H^1_{\FF}(k_\lambda,M)$ under the local Tate pairing $<\,,\,>_\lambda$.
\end{define}

\subsubsection{Selmer structures and Selmer groups} 
\label{sec:selmer structure}
Notation from \S\ref{local conditions} is in effect throughout this section. We will also denote $G_{k_{\lambda}}=\textup{Gal}(\overline{k_{\lambda}}/k_{\lambda})$ by $\mathcal{D}_{\lambda}$, whenever we would like to identify this group by a closed subgroup of $G_{k}=\textup{Gal}(\overline{k}/k)$; namely with a particular decomposition group at $\lambda$ in $G_{k}$. We further define $\mathcal{I}_{\lambda} \subset \mathcal{D}_{\lambda}$ to be the inertia group and $\textup{Fr}_{\lambda} \in \mathcal{D}_{\lambda}/\mathcal{I}_{\lambda}$ to be the arithmetic Frobenius element  at $\lambda$. 

\begin{define}
\label{selmer structure}
A \emph{Selmer structure} $\FF$ on $M$ is a collection of the following data:
\begin{itemize}
\item a finite set $\Sigma(\FF)$ of places of $k$, including all infinite places and primes above $p$, and all primes where $M$ is ramified.
\item for every $\lambda \in \Sigma(\FF)$ a local condition (in the sense of \S\ref{local conditions}) on $M$ (which we view now as a $R[[\mathcal{D}_{\lambda}]]$-module), i.e., a choice of $R$-submodule $$H^1_{\FF}(k_{\lambda},M) \subset H^1(k_{\lambda},M)$$ 
 \end{itemize}
If $\lambda \notin \Sigma(\FF)$ we will also write $H^1_{\FF}(k_{\lambda},M)=H^1_{\textup{f}}(k_{\lambda},M)$, where the module $H^1_{\textup{f}}(k_{\lambda},M)$ is the \emph{finite} part of $H^1(k_{\lambda},M)$, defined as in~\cite{mr02} Definition 1.1.6.
\end{define}

\begin{define}
\label{selmer group}
If $\FF$ is a Selmer structure on $M$, we define the \emph{Selmer module} $H^1_{\FF}(k,M)$ to be the kernel of the sum of the restriction maps $$H^1(\textup{Gal}(k_{\Sigma(\FF)}/k),M) \lra \bigoplus_{\lambda \in \Sigma(\FF)}H^1(k_{\lambda},M)/H^1_{\FF}(k_{\lambda},M)$$ where $k_{\Sigma(\FF)}$ is the maximal extension of $k$ which is unramified outside $\Sigma(\FF)$.
\end{define}

\begin{example}
\label{example:canonical selmer}
In this example we recall~\cite{mr02} Defintion 5.3.2. Let $R=\LL$ be the cyclotomic Iwasawa algebra, and let $\mathbb{M}$ be a free $R$ module endowed with a continuos action of $G_k$, which is unramified outside a finite set of places of $k$.  We define a Selmer structure $\FF_\LL$ on $\mathbb {M}$ by setting $\Sigma(\FF_\LL)$ and $H^1_{\FF_\LL}(k_\lambda, \mathbb {M})=H^1(k_\lambda, \mathbb{M})$ for $\lambda \in \Sigma(\FF_\LL)$. This is what we call the \emph{canonical Selmer structure} on $\mathbb{M}$.

 As in Definition~\ref{def:propagation}, induced Selmer structure on the quotients $\mathbb{M}/I\mathbb{M}$ is still denoted by $\FF_\LL$.  Note that $H^1_{\FF_\LL}(k_\lambda, \mathbb{M}/I\mathbb{M})$ will not usually be the same as $H^1(k_\lambda, \mathbb{M}/I\mathbb{M})$. 
\end{example}
\begin{rem}
\label{rem:canonical selmer}
When $R=\LL$ and $\ \mathbb{M}=T\otimes\LL$ (which is the case of interest in this paper), the Selmer structure $\FFc$ of~\cite{kbb} \S2.1 on the quotients $T\otimes\LL/(f)$ may be identified\footnote{For every prime $\lambda$ of $k$, $H^1_{\FFc}(k_\lambda,T\otimes\LL/(f))$ is the image of $H^1(k_\lambda,T\otimes\LL)$ under the canonical map $H^1(k_\lambda,T\otimes\LL) \ra H^1(k_\lambda,T\otimes\LL/(f))$ by the proofs of~\cite{kbb} Propositions 2.10 and 2.12. $H^1_{\FF_\LL}(k_\lambda,T\otimes\LL/(f))$ is exactly the same thing by its very definition.}, under our hypotheses on $\chi$, with the propagation of $\FF_\LL$  to the quotients $T\otimes\LL/(f)$, for every distinguished polynomial $f$ inside $\LL$.
\end{rem}

\begin{define}
\label{def:selmer triple}
A \emph{Selmer triple} is a triple $(M,\FF,\PP)$ where $\FF$ is a Selmer structure on $M$ and $\PP$ is a set of rational primes, disjoint from $\Sigma(\FF)$.
\end{define}

\begin{rem}
Although one might identify the cohomology groups in our setting (when the Galois module in question is $T\otimes\LL$ with $T=\ZZ_p(1)\otimes\chi^{-1}$, or its quotients by ideals of $\LL$) with certain groups of units, using Kummer theory, we will insist on using the cohomological language for the sake of notational consistency with~\cite{mr02} from which we borrow the main technical results. This way, we also hope that it will be easier to hypothesize our approach for  potential generalizations to other settings.
\subsection{Modifying the local conditions at $p$}
\label{modify-at-p}
In~\cite{kbbstark}, we modify the classical local conditions at primes of $k$ above $p$ to obtain a Selmer structure $\FF_{\mathcal{L}}$ on $T$ (see~\cite{kbbstark} \S1). The objective of this section is to lift the Selmer structure $\FF_{\mathcal{L}}$ to a Selmer structure on $T\otimes\LL$.
\end{rem}

In this section we will make use of the results from the Appendix to determine the structure, as a $\LL$-module, of the semi-local cohomology group $H^1(k_p,T\otimes\LL)$. There should of course be a more direct way (in this setting where $T=\ZZ_p(1)\otimes\chi^{-1}$) to obtain the results on the structure of $H^1(k_p,T\otimes\LL)$, we believe the more general approach via Fontaine's theory of $(\varphi,\Gamma)$-modules might allow our strategy to apply in many other settings.    

Let $k_{\infty}$ denote the cyclotomic $\ZZ_p$-extension of $k$, and $\Gamma=\Gal(k_{\infty}/k)$, as before. Since we assumed  $p$ is unramified in $k/\QQ$, note that $k_{\infty}/k$ is then totally ramified at all primes $\wp$ over $p$. Let $k_{\wp}$ denote the completion of $k$ at $\wp$, and let  $k_{\wp,\infty}$ denote the cyclotomic $\ZZ_p$-extension of $k_{\wp}$. We may identify $\Gal(k_{\wp,\infty}/k_{\wp})$ by $\Gamma$ for all $\wp|p$ and henceforth $\Gamma$ will stand for any of these Galois groups. $\LL=\ZZ_p[[\Gamma]]$ is the cyclotomic Iwasawa algebra, as usual. We also fix a topological generator $\gamma$ of $\Gamma$, and we set $\xx=\gamma-1$ (and occasionally we identify $\LL$ by the power series ring $\ZZ_p[[\xx]]$ in one variable).

Recall that $H^1_{\textup{Iw}}(k_{\wp},T):=\varprojlim_n H^1(k_{\wp,n},T)$, where $k_{\wp,n}$ denotes the unique subfield of $k_{\wp,\infty}$ which has degree $p^n$ over $k_{\wp}$. By Shapiro's Lemma one may canonically identify $H^i_{\textup{Iw}}(k_{\wp},T)$ by $H^i(k_{\wp},T\otimes\LL)$  for all $i \in \ZZ^+$ (\emph{c.f.} ~\cite{colmez-annals} Proposition II.1.1). Let $$H^i_{\textup{Iw}}(k_{p},T):=\oplus_{\wp|p}H^i_{\textup{Iw}}(k_{\wp},T) \hbox{ and } H^i(k_{p},T\otimes\LL):=\oplus_{\wp|p} H^i(k_{\wp},T\otimes\LL)$$
(and these two $\LL$-modules are canonically isomorphic by above argument).

Recall from Appendix \ref{fontaine} that $H_K$ denotes $\Gal(\overline{K}/K_{\infty})$ for any local field $K$. We also set $T_m=T\otimes\LL/(\xx^m)$ and $T_{s,m}=T\otimes\LL/(p^s,\xx^m)$ for $s, m\in\ZZ^+$, following~\cite{kbb} \S2.3.2. 
\begin{prop}
\label{structure-at-p}
\begin{enumerate}
\item[\textbf{(i)}] $H^1(k_{p},T\otimes\LL)=H^1_{\textup{Iw}}(k_{p},T)$ is a free $\LL$-module of rank $r$.
\item[\textbf{(ii)}] The map $H^1(k_{p},T\otimes\LL) \ra H^1(k_{p},T_m)$ is surjective all $m \in \ZZ^+$.
\end{enumerate}
\end{prop}
\begin{proof}
Since we assumed $\chi(\wp)\neq1$ for any $\wp|p$, it follows that $T^{H_{k_{\wp}}}=0$, and thus (i) follows immediately  from Theorem~A.\ref{cohomology-iwasawa-p}.  It also follows that $(T^*)^{G_{k_{\wp}}}=0$ since $\chi(\wp)\neq1$  and hence, by the proof of~\cite{kbb} Lemma~2.11, that $H^2(k_{\wp},T\otimes\LL)=0$. But then $$\textup{coker}\{H^1(k_{p},T\otimes\LL) \ra H^1(k_{p},T_m)\}=\oplus_{\wp|p}H^2(k_{\wp},T\otimes\LL)[\xx^m]=0$$ hence (ii) follows.
\end{proof}

Fix a free $\LL$-direct summand $\LLL$ of $H^1(k_{p},T\otimes\LL)$ which is free of rank one as a $\LL$-module. By Proposition~\ref{structure-at-p} this also fixes a free $\LL/(\xx^m)$-direct summand $\mathcal{L}_m$ of $H^1(k_{p},T_m)$ which is free of rank one (as a $\LL/(\xx^m)$-module). When $m=1$, we denote $\mathcal{L}_1$ simply by $\mathcal{L}$.

\begin{define}
\label{modified selmer structure}
Let $\LLL$ be as above. We define the $\LLL$-\emph{modified Selmer structure} $\FF_{\LLL}$ on $T\otimes\LL$ as
\begin{itemize}
\item $\Sigma(\FF_{\LLL})=\Sigma(\FF_{\LL})$,
\item $H^1_{\FF_{\LLL}}(k_p,T\otimes\LL)= \LLL\subset H^1(k_p,T\otimes\LLL)$,
\item $H^1_{\FF_{\LLL}}(k_{\lambda},T\otimes\LL)=H^1_{\FF_{\LL}}(k_{\lambda},T\otimes\LL)$, for primes $\lambda\nmid p$ of $k$.
\end{itemize}
\end{define}

The induced Selmer structure on the quotients $\mathcal{T}_0:=\{T_{s,m}\}$ will also be denoted by $\FF_{\LLL}$ (except for the induced Selmer structure on $T=T_1=T\otimes\LL/(\xx)$ and its quotients $T/p^sT=T_{s,1}=T\otimes\LL/(p^s,1)$, which we will denote by $\FF_\mathcal{L}$, for notational consistency with~\cite{kbbstark}).
\subsection {Local duality and the dual Selmer structure} We will discuss local duality in full generality. Let $R$ and $M$ be as above, namely $R$ is a complete local noetherian ring and $M$ is a free $R$-module of finite rank which is endowed with a continuous action of $G_k$. Let $M^*=\textup{Hom}(M,\mu_{p^{\infty}})$ be the Cartier dual of $M$. For each prime $\lambda$ of $k$, there is a perfect pairing, called the local Tate pairing $$<,>_{\lambda}: H^1(k_\lambda,M)\times H^1(k_\lambda,M^*)\lra \QQ_p/\ZZ_p$$

Let $\FF$ be a Selmer structure on $M$. For each prime $\lambda$ of $k$, define $H^1_{\FF^*}(k_\lambda,M^*):=H^1_\FF(k_\lambda,M)^\perp$ as the orthogonal complement of $H^1_\FF(k_\lambda,M)$ under the local Tate pairing. The Selmer structure $\FF^*$ on $M^*$ defined in this way will be called the \emph{dual Selmer structure}. 

As in Definition~\ref{selmer group}, dual Selmer structure gives rise to the \emph{dual Selmer group}: $$H^1_{\FF^*}(k,M^*) =\ker\left\{H^1(\textup{Gal}(k_{\Sigma(\FF)}/k),M^*) \lra \bigoplus_{\lambda \in \Sigma(\FF)}\frac{ H^1(k_{\lambda},M^*)}{H^1_{\FF}(k_{\lambda},M^*)}\right\}$$

\subsection{Comparison of Selmer modules}
As our sights are set on Iwasawa's main conjecture over totally real number fields, we first construct the \emph{correct} Iwasawa module: A Selmer group which should relate to the appropriate $p$-adic $L$-function (which has been constructed by Deligne and Ribet~\cite{deligne-ribet} in this setting).

Once this Selmer group is defined, we will use Poitou-Tate global duality to compare it to $H^1_{\FF_{\LLL}^*}(k,(T\otimes\LL)^*)$, dual Selmer group attached to dual $\LLL$-modified Selmer structure. 

\begin{define}
\label{strict Selmer}
\emph{p-strict Selmer structure} $\FF_{\textup{str}}$ on $T\otimes\LL$ is defined by the following data:
\begin{itemize}
\item $\Sigma(\FF_{\textup{str}})=\Sigma(\FF_{\LL})$,
\item $H^1_{\FF_{\textup{str}}}(k_p,T\otimes\LL)=0$,
\item $H^1_{\FF_{\textup{str}}}(k_{\lambda},T\otimes\LL)=H^1_{\FF_{\LL}}(k_{\lambda},T\otimes\LL)=H^1_{\FF_{\LLL}}(k_\lambda,T\otimes\LL)$, for primes $\lambda\nmid p$ of $k$.
\end{itemize}
\end{define}

Hence for the dual Selmer structure $\FF_{\textup{str}}^*$ to $p$-strict Selmer structure $\FF_{\textup{str}}$, we have 
\begin{itemize}
\item $H^1_{\FF_{\textup{str}}^*}(k_p,(T\otimes\LL)^*)=H^1(k_p,(T\otimes\LL)^*)$,
\item $H^1_{\FF_{\textup{str}}^*}(k_{\lambda},(T\otimes\LL)^*)=H^1_{\FF_{\LLL}^*}(k_{\lambda},(T\otimes\LL)^*)$, for primes $\lambda\nmid p$ of $k$.
\end{itemize}

Later in \S\ref{applications}, we explain why $H^1_{\FF_{\textup{str}}^*}(k_p,(T\otimes\LL)^*)$ is the Iwasawa module which (conjecturally\footnote{which already has been proved by Wiles~\cite{wiles-mainconj} using different techniques than ours. Wiles systematically uses Hida's theory of $\Lambda$-adic Hilbert modular forms to construct certain unramified extensions from which he deduces the main conjectures.}) relates to the Deligne-Ribet $p$-adic  $L$-function.

For any $\ZZ_p$-module $A$, let $A^{\vee}:=\textup{Hom}(A,\QQ_p/\ZZ_p)$ denote its Pontryagin dual.  
\begin{prop}\label{exact-seq-stark}
Assume Leopoldt's conjecture. Then the sequence \[
\begin{array}{rl}
  0&\lra H^1_{\FF_{\LLL}}(k,T\otimes\LL)\stackrel{\textup{loc}_p}{\lra} \LLL/\LL \lra   \\
 &H^1_{\FF_{\textup{str}}^*}(k,(T\otimes\LL)^*)^{\vee} \lra (H^1_{\FF_{\LLL}^*}(k,(T\otimes\LL)^*)^{\vee}\lra 0   \\
\end{array}
\]  is exact. 
\end{prop}
\begin{proof}
 First injection follows from Leopoldt's conjecture, and the rest of the sequence is exact by Poitou-Tate global duality (which is a statement of class field theory in our setting) See also~\cite{r00} Theorem I.7.3, proof of Theorem III.2.10 and \cite{deshalit} \S III.1.7.
\end{proof}

\subsection{Kolyvagin systems for the $\LLL$-modified Selmer triple - I}
\label{kolsys-1}
We first recall our notation. $\chi$ is an $even$ character of Gal$(\overline{k}/k)$ ($i.e.$ it is trivial on all complex conjugations inside Gal$(\overline{k}/k)$) into $\ZZ_p^{\times}$ that has finite order, and $L$ is the fixed field of $\hbox{ker}(\chi)$, inside $\overline{k}$. $f_{\chi}$ is the conductor of $\chi$, and $\Delta$ is the Galois group $\hbox{Gal}(L/k)$ of the extension $L/k$. We assume that $p\nmid f_{\chi}$ and $\chi(\wp)\neq1$ for any prime $\wp$ of $k$ over $p$. We assume further that $p$ is unramified in $k/\QQ$.  $T$ is the Gal$(\overline{k}/k)$-representation $\ZZ_p(1)\otimes\chi^{-1}$. 
\begin{rem}
\label{hypo holds}
$T$ clearly satisfies the hypotheses \hone-\hfour\, ($\QQ$ replaced by $k$) of~\cite{kbb} \S2.2. Observe also that the below versions of the hypotheses \htam\, and \hsez\, will hold for $T$ as well:
\begin{enumerate}
\item[(\htamF)] $A ^{I_{F_{\lambda}}}$ is divisible for every prime $\lambda \nmid p$ of $k$.
\item[(\hsezF)]  $(T^*)^{G_{F_{\wp}}} =0$ for primes $\wp\mid p$ of $k$.
\end{enumerate}
\end{rem}

Recall also the definition of the collection $\mathcal{T}_0=\{T_{s,m}\}$. By definition, local conditions on $T_{s,m}$ at primes $\lambda\nmid p$ of $k$ determined by $\FF_{\LLL}$ will coincide with the local conditions determined by $\FF_\LL$, and therefore, by the proofs of~\cite{kbb} Propositions 2.10 and 2.12, they will also coincide with the local conditions determined by $\FFc$ since  \htamF\, and \hsezF\, hold true.  

\begin{prop}
\label{prop:cart-L}
\begin{enumerate}
\item[\textbf{(i)}] The Selmer structure $\FF_{\LLL}$ is cartesian on $\mathcal{T}_0$, in the sense of \cite{kbb} Definition 2.4.
\item[\textbf{(ii)}] The core Selmer rank (c.f.~\cite{mr02} Definition 4.1.11 for a definition) $\XX(T,\FF_\mathcal{L})$ of the Selmer structure on $T$ induced from $\FF_{\LLL}$ is one.
\end{enumerate}
\end{prop}
\begin{proof}
(ii) is Proposition 1.8 of~\cite{kbbstark}. 

Since $\FFc$ and $\FF_{\LLL}$ determine the same local conditions at places $v\nmid p$, the local conditions are cartesian at $v\nmid p$ by~\cite{kbb} Proposition 2.10. Therefore, it suffices to check that $\FF_{\LLL}$ is cartesian on $\mathcal{T}_0$ at $p$, i.e. we need to verify properties $\mathbf{C.1}$-$\mathbf{C.3}$ of~\cite{kbb} Definition 2.4 for the local conditions determined by $\FF_{\LLL}$ on $\mathcal{T}_0$.

Property $\mathbf{C.1}$  holds by definition of $\FF_{\LLL}$ on $\mathcal{T}_0$, and property $\mathbf{C.3}$ follows easily from~\cite{mr02} Lemma 3.7.1 (which applies since $\mathcal{L}_m$ is a direct summand of $H^1(k_p,T_m)$, i.e. $H^1(k_p,T_m)/\Ll_m$ is a free $\LL/(\xx^m)$-module for all $m$).

We now verify $\mathbf{C.2}$. Let $\Ll_{s,m}$ be the image of $\Ll_m$ under the reduction map $$H^1(k_p,T_m) \lra H^1(k_p,T_{s,m})$$ for $s\in\ZZ^+$. It is easy to see that $\Ll_{s,m}$ (resp. $H^1(k_p,T_{s,m})/\Ll_{s,m}$) is a free $\LL/(p^s,\xx^m)$-module of rank one (resp. of rank $[k:\QQ]-1$). We need to check that the map $$H^1(k_p,T_{s,m})/\Ll_{s,m} \stackrel{[\xx^{M-m}]}{\lra} H^1(k_p,T_{s,M})/\Ll_{s,M} $$ induced from the map $[\xx^{M-m}]:\LL/(\xx^m)\ra\LL/(\xx^M)$ is injective for all $M\geq m$. But this is evident since $H^1(k_p,T_{s,m})/\Ll_{s,m}$ (resp. $H^1(k_p,T_{s,M})/\Ll_{s,M}$) is a free $\LL/(p^s,\xx^m)$-module (resp. a free $\LL/(p^s,\xx^M)$-module) of rank $[k:\QQ]-1$.
\end{proof}

 As a corollary of Proposition~\ref{prop:cart-L} (and the proof of~\cite{kbb} Theorem 3.23) we obtain
\begin{thm}
\label{main-stark}
The $\LL$-module of Kolyvagin Systems $\overline{\KS}(T\otimes\LL,\FF_{\LLL})$ for the Selmer structure $\FF_{\LLL}$ on $T\otimes\LL$  is free of rank one. Further, the map $$\overline{\KS}(T\otimes\LL,\FF_{\LLL}) \lra \overline{\KS}(T,\FF_{\mathcal{L}})$$ is surjective.
\end{thm}
See~\cite{kbb} \S\S3.1-3.2 for a precise definition of the module of Kolyvagin systems $\overline{\KS}(T\otimes\LL,\FF_{\LLL})$.
\begin{proof}
This is exactly~\cite{kbb} Theorem 3.23 for $T=\ZZ_p(1)\otimes\chi^{-1}$ and the base field is $k$ instead of $\QQ$. (The proof of \cite{kbb} Theorem 3.23 generalizes word by word to the case when the base field is different from $\QQ$. The only technical point to verify is Proposition~\ref{prop:cart-L} above.)
\end{proof}

In the next section, we will show how to obtain a generator of the cyclic $\LL$-module $\overline{\KS}(T\otimes\LL,\FF_{\LLL})$ using Stark elements of Rubin.
\newpage
\section{Kolyvagin systems of Stark units} 
\label{kolsys-starkunits}
In this section we review Rubin's integral refinement of Stark conjectures and construct Kolyvagin systems for the modified Selmer structure $(T\otimes\LL,\FF_{\LLL})$ (whose existence was proved unconditionally in the previous section, building on ideas from~\cite{kbb}) using Stark elements of Rubin.

For the rest of this paper we assume this refined version of Stark conjectures.

Before we give an outline of Rubin's conjectures, we set some notation. Assume $k,k_\infty, \chi,f_\chi$ and $L$ are as above. For a cycle $\tau$ of the number field $k$, let $k(\tau)$ be the maximal $p$-extension inside the ray class field of $k$ modulo $\tau$. Let also $k_n$ be the unique subextension of $k_\infty/k$ which has degree $p^n$ over $k$. For any other number field $F$, we define $F(\tau)$ as the composite of $k(\tau)$  and $F$; and $F_n$ as the composite of $k_n$ and $F$. Let $$\mathcal{K} =\{L_n(\tau): \tau \hbox{ is a (finite) cycle of } k \hbox{ prime to } f_{\chi}p\}$$ and  $$\mathcal{K}_0 =\{k_n(\tau): \tau \hbox{ is a (finite) cycle of } k \hbox{ prime to } f_{\chi}p\}$$ be collections of abelian extensions  of $k$. 

\subsection{Stark elements and Euler systems (of rank $r$) for $\ZZ_p(1)$}
\label{sec:stark-euler}
Fix a finite set $S$ of places of $k$ that does \emph{not} contain $p$, but contains all infinite places $S_{\infty}$, all places $\lambda$ that divide the conductor $f_{\chi}$ of $\chi$. Assume that $|S| \geq r+1$. For each $K \in \kk$ let $S_K=S \cup \{\hbox{places of } k \hbox{ at which }K \hbox{ is ramified}\}$ be another set of places of $k$. Let $\mathcal{O}_{K,S_K}^{\times}$ denote the $S_K$ units of $K$ (i.e. elements of $K^\times$ that are units away from the primes of $K$ above $S_K$), and $\Delta_K$ (\emph{resp.} $\delta_K$) denote $\hbox{Gal}(K/k)$ (\emph{resp.} $|\hbox{Gal}(K/k)|$).  Conjecture $B^{\prime}$ of~\cite{ru96} predicts the existence of certain elements $$\varepsilon_{K,S_K} \in \Lambda_{K,S_K} \subset \frac{1}{\delta_K}{\bigwedge^r} \mathcal{O}_{K,S_K}^{\times}$$ where $\Lambda_{K,S_K}$ is defined in \S 2.1 of~\cite{ru96} and has the property that for any homomorphism $$\tilde{\psi} \in \hbox{Hom}_{\QQ_p[\Delta_K]}(\bigwedge^{r}\mathcal{O}_{K,S_K}^{\times,\wedge} \otimes \QQ_p,\mathcal{O}_{K,S_K}^{\times, \wedge} \otimes \QQ_p)$$ that is induced from a homomorphism $$\psi \in \hbox{Hom}_{\ZZ_p[\Delta_K]}(\bigwedge^{r}\mathcal{O}_{K,S_K}^{\times,\wedge}, \mathcal{O}_{K,S_K}^{\times, \wedge})$$ we have $\psi(\Lambda_{K,S_K}) \subset \mathcal{O}_{K,S_K}^{\times, \wedge}$. 
\begin{rem}\label{rem:T}
In fact, Rubin's conjecture predicts that these elements should be inside  $\frac{1}{\delta_K}{\bigwedge^r} \mathcal{O}_{K,S_K,T}^{\times}$ where $T$ is a finite set of primes disjoint from $S_K$ chosen so that the group $\mathcal{O}_{K,S_K,T}^{\times}$ of units which are congruent to 1 modulo all the primes in $T$ is torsion-free. However, in our case any set $T$ which contains a prime other than 2 will suffice (since all the fields that appear in our paper are totally real). 

Further, $T=\{\qq\}$ may be chosen in a way that the extra factors that appear in the definition of the modified $zeta\,function$ for $K$ ($c.f.$ \cite{ru96}~\S1 for more detail on these zeta functions) will be prime to $p$, when they are evaluated at $0$. (This could be accomplished, for example, by choosing $\qq$ so that $\mathbf{N}\qq-1$ is prime to $p$.) We note that for such a chosen $T$, we have $\mathcal{O}_{K,S_K,T}^{\times,\wedge}=\mathcal{O}_{K,S_K}^{\times,\wedge}$ , for example by the exact sequence (1) in~\cite{ru96}. Since in our paper we will work with the $p$-adic completion of the units, we will safely exclude $T$ from our notation.

One minor point arises because of the appearance of the set of primes $T=\{\qq\}$: One should remove fields $k_n(\tau)$ (and $L_n(\tau)$) from the collections $\mathcal{K}_0$ (and $\mathcal{K}$, respectively) for which $\qq|\tau$. This also is not a problem at all for our purposes.
\end{rem} 
The elements $\varepsilon_{K,S_K}$ (which we call Stark elements) satisfy the distribution relation to be satisfied by an \emph{Euler system of rank r} (in the sense of~\cite{pr-es}) as $K$ runs through the extensions $\mathcal{K}$, by Proposition 6.1 of~\cite{ru96}. We denote the image of $\varepsilon_{K,S_K}$ inside the $\ZZ_p$-module $\frac{1}{\delta_K}\bigwedge^r \mathcal{O}_{K,S_K}^{\times,\wedge}$ also by $\varepsilon_{K,S_K}$. Since $S$ is fixed (therefore $S_K$, too), we will often drop $S$ or $S_K$ from notation and denote $\varepsilon_{K,S_K}$ by $\varepsilon_{K}$, or sometimes use $S$ instead of $S_K$ and denote $\mathcal{O}_{K,S_K}$ by $\mathcal{O}_{K,S}$. 

For any number field $K$ one may identify $H^1(K,\ZZ_p(1))$ with $$K^{\times,\wedge} := \lim_{\stackrel{\longleftarrow}{n}}K^{\times}/(K^{\times})^{p^n}$$ by Kummer theory. Under this identification we view  $\varepsilon_{K,S_K}$ as an element of $\bigwedge^r H^1(K,\ZZ_p(1))$. Distribution relation satisfied by the Stark elements (\cite{ru96} Proposition 6.1) shows that the collection $\{\varepsilon_{K,S_K}\}_{K\in\mathcal{K}}$ is an Euler system of rank $r$ in the sense of~\cite{pr-es}. Following the formalism of~\cite{r00} \S II.4, we \emph{twist} this Euler system to obtain an Euler system for the representation $T=\ZZ_p(1)\otimes\chi^{-1}$.
\subsection{Twisting by the character $\chi$}
\label{sec:twisting}
 Set $\Gamma_n:=\Gal(k_n/k)$, $G_{\tau}:=\Gal(k(\tau)/k)$,  $\Delta_{\tau}:=\Gal(L(\tau)/k)=G_\tau\times\Delta$, $G^{(n)}_{\tau}:=\Gal(k_n(\tau)/k)=G_\tau\times\Gamma_n$, which is the $p$-part of $\Delta^{(n)}_{\tau}:=\hbox{Gal}(L_n(\tau)/k)\cong G^{(n)}_{\tau} \times \Delta=G_\tau\times\Gamma_n\times\Delta$. (These canonical factorizations of the Galois groups follow easily from the fact that  $|\Delta|$ is prime to $p$ and from ramification considerations.)
 
 Let $\chi$ be as above. Let $\eee_{\chi}$ denote the idempotent $\frac{1}{|\Delta|}\sum_{\delta \in \Delta}\chi(\delta)\delta^{-1}$. We regard this element as an element of the groups rings $\ZZ_p[\Delta^{(n)}_{\tau}]$ via above factorizations.
 
 For any cycle $\tau$ which is prime to $pf_{\chi}$ we define 
\begin{align} \label{def:twist}
\varepsilon_{L_n(\tau)}^{\chi}:=\epsilon_{\chi}\varepsilon_{L_n(\tau),S} \in  \epsilon_{\chi}\bigwedge^r H^1(L_n(\tau),\ZZ_p(1))&=\\ \label{eqn*}
= \bigwedge^r \eee_\chi H^1(L_n(\tau),\ZZ_p(1))&=\\
= \bigwedge^r H^1(L_n(\tau),\ZZ_p(1))^{\chi}&
\end{align}                          
We also note that (\ref{eqn*}) above 
\begin{align*}
\left(\bigwedge^r H^1(L_n(\tau),\ZZ_p(1))\right)^{\chi}=\epsilon_{\chi}\bigwedge^r H^1(L_n(\tau),\ZZ_p(1)) &=\\ =\bigwedge^r \eee_\chi H^1(L_n(\tau),\ZZ_p(1))&=\\
=\bigwedge^r H^1(L_n(\tau),\ZZ_p(1))^{\chi} &
\end{align*}
holds simply because $\epsilon^{r}_{\chi}=\epsilon_{\chi}$.
 
 Now Hochschild-Serre spectral sequence gives rise to an exact sequence 
\begin{equation}
\label{eq:twist}
H^1(\Delta,T) \lra H^1(k_n(\tau),T)\lra H^1(L_n(\tau),T)^{\Delta}\lra H^2(\Delta,T)
\end{equation}
where $H^1(L_n(\tau),T)^{\Delta}$ stands for the largest submodule of $H^1(L_n(\tau),T)$ on which $\Delta$ acts trivially. On the other hand, since $|\Delta|$ is prime to $p$, it follows that very first and the very last terms in~(\ref{eq:twist}) vanish. Therefore, we have an isomorphism $$H^1(k_n(\tau),\ZZ_p(1)\otimes\chi^{-1})\lra H^1(L_n(\tau),\ZZ_p(1)\otimes\chi^{-1})^{\Delta}$$
On the other hand, since $G_{L_n(\tau)}$ is in the kernel of $\chi$ $$H^1(L_n(\tau),\ZZ_p\otimes\chi^{-1}) \cong H^1(L_n(\tau),\ZZ_p(1))\otimes\chi^{-1}$$ hence 
\begin{equation}
\label{twist-isom}
H^1(k_n(\tau),T) \stackrel{\sim}{\lra} H^1(L_n(\tau),T)^{\Delta} \cong H^1(L_n(\tau),\ZZ_p(1))^{\chi} 
\end{equation}
where, as we have set earlier, $H^1(L_n(\tau),\ZZ_p(1))^{\chi}$ is the $\chi$-isotypic part of $H^1(L_n(\tau),\ZZ_p(1))$. This induces an isomorphism
\begin{equation}\label{eq:main twist}
\bigwedge^r H^1(k_n(\tau),T) \stackrel{\sim}{\lra} \bigwedge^r H^1(L_n(\tau),\ZZ_p(1))^{\chi} 
\end{equation}
The inverse image of the element $\varepsilon_{L_n(\tau)}^{\chi}$  (which we defined in~(\ref{def:twist})) under the isomorphism~(\ref{eq:main twist}) above will be denoted by $\varepsilon_{k_n(\tau)}^{\chi}$. The collection $\{\varepsilon_{K}^{\chi}\}_{K\in\mathcal{K}_0}$ we be called the \emph{Stark element Euler system of rank r}.

Next we construct an \emph{Euler system of rank one} (i.e. an Euler system in the sense of~\cite{r00}) using ideas from~\cite{ru96, pr-es} (which are originally due to Rubin). Main point is that, if one applies the arguments of~\cite{ru96,pr-es} directly, all one could get at the end (after applying Kolyvagin's decent) is a $\LL$-adic Kolyvagin system for the much coarser Selmer structure $\FFc$ on $T\otimes\LL$. In Section \S\ref{kolsys-2}, we overcome this difficulty and obtain a $\LL$-adic Kolyvagin system  for the refined Selmer structure $\FF_{\LLL}$ on $T\otimes\LL$.

\subsection{Choosing homomorphisms}
\label{sec:homs}
For any field $K\in\kk_0$, we define $\Delta_K:=\Gal(K/k)$. Using elements of $$\varprojlim_{K \in \kk_0} \bigwedge^{r-1}\hbox{Hom}_{\ZZ_p[\Delta_K]}(H^1(K,T), \ZZ_p[\Delta_K])$$  (or the images of the elements of $$\varprojlim_{K \in \kk_0} \bigwedge^{r-1}\hbox{Hom}_{\ZZ_p[\Delta_K]}(H^1(K_p,T), \ZZ_p[\Delta_K])$$ under the canonical map 
$$\xymatrix{
\varprojlim_{K \in \kk_0} \bigwedge^{r-1}\hbox{Hom}_{\ZZ_p[\Delta_K]}(H^1(K_p,T), \ZZ_p[\Delta_K]) \ar[d]\\
\varprojlim_{K \in \kk_0} \bigwedge^{r-1}\hbox{Hom}_{\ZZ_p[\Delta_K]}(H^1(K,T), \ZZ_p[\Delta_K])
}$$   induced from localization at $p$) and the Stark elements above, Rubin~\cite{ru96} \S6 (see also~\cite{pr-es} for an application of this idea in a more general setting) shows how to obtain Euler system (in the sense of~\cite{r00}) for the representation $T$ of $G_k$. In this section, we show how to choose these homomorphisms carefully so that the resulting Euler system gives rise to a Kolyvagin system for the $\LLL$-modified Selmer structure $\FF_{\LLL}$ on $T\otimes\LL$.

Recall that one may identify  $H^1(k_n(\tau),T)$ with $(L_n(\tau)^{\times})^{\chi}$, using (\ref{twist-isom}) and Kummer theory. Similarly, one may identify the semi-local cohomology group $H^1(k_n(\tau)_p,T)$ with $(L_n(\tau)_p^{\times})^{\chi}$, where $L_n(\tau)_p:=L_n(\tau)\otimes\QQ_p$. Let $V_{L_n(\tau)}$ denote the $p$-adic completion of the local units of $L_n(\tau)_p$.

On the other hand, the proof of~\cite{r00} Proposition III.2.6(ii) shows that (since we assume $\chi(\textup{Frob}_\wp)\neq1$ for any prime $\wp|p$ of $k$):
\begin{lemma}
\label{at p: relaxed=unr}
$$H^1(k_n(\tau)_p,T)\cong(L_n(\tau)_p^{\times})^{\chi}\cong V_{L_n(\tau)}^{\chi}$$ for all $k_n(\tau) \in \kk_0$.
\end{lemma}

We remark that all fields that appear in this paper (namely elements of the collection $\kk$) are all totally real. Further, since we assumed that $p$ is prime to $f_\chi$ it follows that $L(\tau)/k$ is unramified at all primes above $p$. Therefore Krasner's Lemma~\cite{krasner} on the structure of $1$-units implies:
\begin{lemma}\textup{(Krasner)}
\label{krasner-1}
$V_{L(\tau)}$ is a free $\ZZ_p[[\Delta_\tau]]$-module of rank $r$ \textup{(}where $r=[k:\QQ]$\textup{)}.
\end{lemma}
\begin{cor}
\label{krasner}
$H^1(k(\tau)_p,T)=V_{L(\tau)}^{\chi}$ is a free $\ZZ_p[[G_\tau]]$-module of rank $r$.
\end{cor}
On the other hand, Theorem A.\ref{cohomology-iwasawa-p} shows that (which applies to determine the structure of semi-local cohomology as in \S\ref{modify-at-p} and Proposition~\ref{structure-at-p}(ii), since $k(\tau)$ is still unramified at all primes above $p$. We also remark that $T^{H_{k(\tau)_{\mathcal{Q}}}}=0$ as well, by the proof of~\cite{r00} Lemma 4.2.5(i) (see also the remark following~\cite{r00} Conjecture 8.2.6), for every $k(\tau) \in \kk_0$ and every prime $\mathcal{Q}$ of $k(\tau)$ above $p$; and $H_{k(\tau)_{\mathcal{Q}}}$ is defined at the beginning of Appendix~\ref{fontaine}):
\begin{prop}
\label{Lambda-rank}
\begin{enumerate}
\item[(i)] $H^1(k(\tau)_p,T\otimes\LL)$ is a free $\LL$-module of rank $[k(\tau):\QQ]=r\cdot |G_\tau|$.
\item[(ii)] The canonical map  $H^1(k(\tau)_p,T\otimes\LL)\lra  H^1(k(\tau)_p,T)$ is surjective.
\end{enumerate}
\end{prop}
\begin{rem}
\label{rem:surjective}
The statement of Proposition~\ref{Lambda-rank}(ii) is equivalent to saying that the natural map $$H^1(k(\tau)_p,T\otimes\LL)/(\gamma-1)\cdot H^1(k(\tau)_p,T\otimes\LL)\lra  H^1(k(\tau)_p,T)$$ is an isomorphism.
\end{rem}
Let $\LL_\tau:=\ZZ_p[[G_\tau\times\LL]]=\LL\otimes_{\ZZ_p}\ZZ_p[G_\tau]$ and let $\mathbb{M}$ be the composite of all the fields inside $\kk_0$.
\begin{prop}
\label{big-free}
The $\LL_\tau$-module $H^1(k(\tau)_p,T\otimes\LL)$ is free of rank $r$.
\end{prop}
\begin{proof}
Let $e=\{e_1,\dots,e_r\}$ be a $\ZZ_p[G_\tau]$-basis for $H^1(k(\tau)_p,T)$. By Nakayama's lemma, there exists a set of generators $\mathbb{E}=\{\mathbb{E}_1,\dots,\mathbb{E}_r\}$ of the $\LL_\tau$-module $H^1(k(\tau)_p,T\otimes\LL)$, which lifts  the basis $e$ (under the surjection $$H^1(k(\tau)_p,T\otimes\LL)\lra  H^1(k(\tau)_p,T)$$ of Proposition~\ref{Lambda-rank} and Remark~\ref{rem:surjective}).

Assume that there is a relation \begin{equation}\label{relation}\Sigma_{i=1}^r  a_i\mathbb{E}_i=0\end{equation} with $a_i\in\LL_\tau$. Since $H^1(k(\tau)_p,T\otimes\LL)$ is $\LL$-torsion free (by Proposition~\ref{Lambda-rank}), we may assume without loss of generality that there is a $j\in\{1,\dots,r\}$ such that $a_j\notin (\gamma-1)$. This means, however, relation~(\ref{relation}) reduced modulo $(\gamma-1)$ gives rise to a non-trivial relation among $\{e_1,\dots,e_r\}$ over $\ZZ_p[G_\tau]$. This is a contradiction, which proves that the set of generators $\mathbb{E}$ is a $\LL_\tau$-basis for $H^1(k(\tau)_p,T\otimes\LL)$.
\end{proof}
\begin{cor}
\label{huge-free}
The $\LL\otimes_{\ZZ_p}\ZZ_p[[\Gal(\mathbb{M}/k)]]$-module $$\mathbb{V}=\varprojlim_{\tau} H^1(k(\tau)_p,T\otimes\LL)$$ is free of rank $r$.
\end{cor}
\begin{proof}
Immediate after Proposition~\ref{big-free}.
\end{proof}
\begin{cor}
\label{cor:free}
The $\ZZ_p[G_\tau^{(n)}]$-module $H^1(k_n(\tau)_p,T)$ is free of rank $r$.
\end{cor}
\begin{proof}
Follows from Proposition~\ref{big-free} and by the fact that the map $H^1(k(\tau)_p,T\otimes\LL)=H^1_{\textup{Iw}}(k(\tau)_p,T)\lra H^1(k_n(\tau)_p,T)$ is surjective (as the relevant $H^2$ vanishes).
\end{proof}
 
Choose a $\LL\otimes_{\ZZ_p}\ZZ_p[[\textup{Gal}(M/k)]]$-line $\mathbb{L}$  inside $\vv$ such that the quotient $\vv/\mathbb{L}$ is also a free $\LL\otimes_{\ZZ_p}\ZZ_p[[\textup{Gal}(M/k)]]$-module (of rank $r-1$).
\begin{define}
\label{lines}
For all $k_n(\tau)=K \in \kk_0$ let $\al_{K}$ be the image of $\mathbb{L}$ under the (surjective) projection map \, $\vv \ra H^1(K_p,T)$. 
\end{define}

Note that $\al_K$ are free $\ZZ_p[\textup{Gal}(K/k)]$-modules of rank \emph{one} for all $K \in \kk$, and that $(\al_K)^{\textup{Gal}(K/K^{\prime})}=\al_{K^{\prime}}$ for all $K^{\prime} \subset K$.  When $\tau=1$ (i.e. when $K=k$), denote $\al_K$ by only $\al$.

We will often write  $\al_{n}^{\tau}$ for $\al_{k_n(\tau)}$. We also denote the image of $\all$ under the projection $\vv\ra H^1(k(\tau)_p,T\otimes\LL)$ by $\all_\infty^\tau$. When $\tau=1$, we simply write $\LLL$ for $\all_\infty^\tau$, and $\al_n$ instead of $\al_n^\tau$.

We define 
\begin{align*}
\bigwedge^{r-1}\textbf{Hom}(\mathbb{V},&\,\LL\otimes_{\ZZ_p}\ZZ_p[[\textup{Gal}(M/k)]]):=\\ & \lim_{\stackrel{\longleftarrow}{K \in \kk_0}} \bigwedge^{r-1}\hbox{Hom}_{\ZZ_p[\Delta_K]}(H^1(K_p,T), \ZZ_p[\Delta_K]) \end{align*} where the inverse limit is with respect to the natural maps induced from the inclusion map $H^1(K_p,T) \hookrightarrow H^1(K^{\prime}_p,T)^{\hbox{\tiny Gal}(K^{\prime}/K)}$ (this injection could be easily checked after identifying these cohomology groups with groups of local units) and the isomorphism
 \begin{align*}
  \ZZ_p[\Delta_{K^{\prime}}]^{\hbox{\tiny Gal}(K^{\prime}/K)} & \tilde{\lra}\ZZ_p[\Delta_K]\\
\mathbf{N}^{K^{\prime}}_K &\longmapsto 1
 \end{align*} for $K \subset K^{\prime}$.
 
Localization at $p$ gives rise to a map $H^1(K,T)\stackrel{\textup{loc}_p}{\lra}H^1(K_p,T)$, which  induces a canonical  map 
\begin{align*}
\bigwedge^{r-1}\textbf{Hom}(\mathbb{V},\,\LL\otimes_{\ZZ_p}&\ZZ_p[[\textup{Gal}(M/k)]])\lra\\& \varprojlim_{K\in\kk_0} \bigwedge^{r-1}\textup{Hom}_{\ZZ_p[\Delta_K]}(H^1(K,T), \ZZ_p[\Delta_K])
\end{align*}

The image of $\Phi \in \hhh$ under this map will still be denoted by $\Phi$. 
\begin{prop}
\label{def: es}
Let $$\{\phi_{K}\}_{_{K\in\kk_0}}=\Phi \in \hhh$$ and let $$H^1(K,T) \ni \varepsilon_{K,\Phi}^{\chi}:=\phi_K(\varepsilon_{K}^{\chi}),$$ where $\{\varepsilon_{K}^{\chi}\}_{_{K\in\kk_0}}$ is the Stark element Euler system of rank $r$. Then the collection $\{\varepsilon_{K,\Phi}^{\chi}\}_{_{K\in\kk_0}}$ is an Euler system for the $G_k$-representation $T=\ZZ_p(1)\otimes\chi^{-1}$ in the sense~\cite{r00} Definition 2.1.1. The collection $\{\varepsilon_{K,\Phi}^{\chi}\}_{_{K\in\kk_0}}$ will be called the Euler system of $\Phi$-Stark elements.
\end{prop}

We will sometimes write $\{\varepsilon_{k_n(\tau),\Phi}^{\chi}\}_{_{n,\tau}}$ for the Euler system $\{\varepsilon_{K,\Phi}^{\chi}\}_{_{K\in\kk_0}}$.
\begin{proof}
This is~\cite{ru96} Proposition 6.6. See also~\cite{pr-es} for a more general treatment of this idea of Rubin (which in general shows how to obtain an Euler system, in the sense of~\cite{r00}, starting from an Euler system of rank $r$, in the sense of~\cite{pr-es}).
\end{proof}
 \begin{prop}
 \label{krasner-2}
 For any $K \in \kk_0$ the projection map \begin{align*}\hhh &\lra\\
 & \bigwedge^{r-1}\textup{Hom}_{\ZZ_p[\Delta_K]}(H^1(K_p,T), \ZZ_p[\Delta_K])\end{align*} is surjective.
\end{prop} 
\begin{proof}
Obvious after Corollary~\ref{cor:free}.
\end{proof}

If one applies the \emph{Euler system to Kolyvagin system map} of Mazur and  Rubin (\emph{c.f } \cite{mr02} Theorem 5.3.3) on the Euler system of $\Phi$-Stark elements, all one gets is a $\LL$-adic Kolyvagin system for the coarser Selmer structure $\FFc$ on $T\otimes\LL$. However, below we will choose these homomorphisms $\Phi$ carefully so that the resulting Kolyvagin system will in fact be a Kolyvagin system for the $\LLL$-modified Selmer structure.

\begin{define}
\label{hli}
We say that an element $$\{\phi_n^{\tau}\}_{_{n,\tau}}=\Phi \in \hhh$$ satisfies $\hli$ if for any $K=k_n(\tau)\in \kk_0$ one has $\phi_n^{\tau}(\bigwedge^{r} H^1(K_p,T)) \subset \al_n^{\tau}$.
\end{define}

We will next construct a specific element $$\Phi_0^{(\infty)} \in \hhh$$ that satisfies $\hli$, and such that $\Phi_0^{(\infty)}$ lifts the element $\Phi_0$ of~\cite{kbbstark} \S2.3, under the (surjective) map
\begin{align*}
\bigwedge^{r-1}\textbf{Hom}(\mathbb{V},&\,\LL\otimes_{\ZZ_p}\ZZ_p[[\textup{Gal}(M/k)]])\lra\\ &\varprojlim_\tau \bigwedge^{r-1}\hbox{Hom}_{\ZZ_p[G_\tau]}(H^1(k(\tau)_p,T), \ZZ_p[G_\tau]) \end{align*} 
which was used in~\cite{kbbstark} to construct a \emph{primitive} Kolyvagin system for the Selmer structure $\FF_{\mathcal{L}}$ on $T$.

Fix a basis $$\{\Psi_{\mathbb{L}}^{(i)}\}_{i=1, \dots,r-1}$$ of the free (of rank $r-1$) $\LL\otimes_{\ZZ_p}\ZZ_p[[\textup{Gal}(M/k)]]$-module $$\textup{Hom}_{\LL\otimes_{\ZZ_p}\ZZ_p[[\textup{Gal}(M/k)]]}(\vv/\all,\LL\otimes_{\ZZ_p}\ZZ_p[[\textup{Gal}(M/k)]])$$ This then fixes  a basis $\{\psi_{\al_n^{\tau}}^{(i)}\}_{i=1}^{r-1}$ for the free (of rank $r-1$) $\ZZ_p[G_{\tau}\times\Gamma_n]$-module $\textup{Hom}_{\ZZ_p[G_{\tau}\times\Gamma_n]}\left(H^1(k_n(\tau)_p,T)/\al_n^{\tau}, \,\ZZ_p[G_{\tau}\times\Gamma_n]\right)$ for all $k_n(\tau) \in \kk_0$; such that $\{\psi_{\al_n^{\tau}}^{(i)}\}_{\tau}$ are compatible with respect to the surjections $$\xymatrix{\textup{Hom}_{\ZZ_p[\Delta_K]}(H^1(K_p,T)/\al_K,\ZZ_p[\Delta_K]) \ar[d]\\ \textup{Hom}_{\ZZ_p[\Delta_{K^{\prime}}]}(H^1(K^{\prime}_p,T)/\al_{K^{\prime}},\ZZ_p[\Delta_{K^{\prime}}])}$$ for all $K^\prime=k_{n^\prime}(\tau^\prime),K=k_n(\tau) \in \kk_0$ with $K^\prime\subset K$. Note that the homomorphism $$\bigoplus_{i=1}^{r-1}\psi_{\al_n^{\tau}}^{(i)}:H^1(k_n(\tau)_p,T)/\al_n^{\tau} \lra \ZZ_p[G_{\tau}\times\Gamma_n]^{r-1}$$ is an isomorphism of $\ZZ_p[G_{\tau}\times\Gamma_n]$-modules, for all $\tau$ and $n$.

Let $\psi_{n, \tau}^{(i)}$ denote the image of $\psi_{\al_n^{\tau}}^{(i)}$ under the canonical injection $$\xymatrix{\textup{Hom}_{\ZZ_p[G_{\tau}\times\Gamma_n]}(H^1(k_n(\tau)_p,T)/\al_n^{\tau},\,\ZZ_p[G_{\tau}\times\Gamma_n]) \ar @{^{(}->}[d]\\ \textup{Hom}_{\ZZ_p[G_{\tau}]}(H^1(k_n(\tau)_p,T),\ZZ_p[G_{\tau}\times\Gamma_n])}$$ Note then that $$\Psi_n^{\tau}:=\bigoplus^{r-1}_{i=1}\psi_{n,\tau}^{(i)}: H^1(k_n(\tau)_p,T) \lra \ZZ_p[G_\tau\times\Gamma_n]^{r-1}$$ is surjective and $\textup{ker}(\Psi_n^{\tau})=\al_n^{\tau}$.


Define $$\phi_n^{\tau} := \psi_{n,\tau}^{(1)}\wedge\psi_{n,\tau}^{(2)} \wedge\dots \wedge\psi_{n,\tau}^{(r-1)} \in \bigwedge^{r-1}\textup{Hom}(H^1(k_n(\tau)_p,\ZZ_p[G_{\tau}\times\Gamma_n]) $$ (When $\tau=1$, we drop $\tau$ from the notation and simply write $\phi_n$ for $\phi_n^\tau$, etc.)  Note once again that for $\tau^{\prime}|\tau$ and $n^\prime\leq n$, $\phi_n^{\tau}$ maps to $\phi_{n^\prime}^{\tau^{\prime}}$ under the surjective (by Proposition~\ref{krasner-2} and its proof) homomorphism $$\xymatrix{\bigwedge^{r-1}\textup{Hom}(H^1(k_n(\tau)_p,T),\ZZ_p[G_{\tau}\times\Gamma_n]) \ar@{->>}[d]\\ \bigwedge^{r-1}\textup{Hom}(H^1(k_{n^\prime}(\tau^{\prime})_p,T),\ZZ_p[G_{\tau^\prime}\times\Gamma_{n^\prime}])}$$ Therefore $\Phi_0^{(\infty)}:=\{\phi_n^{\tau}\}_{_{n,\tau}}$ may be regarded as an element of the group $\hhh$.
\begin{prop}
\label{homs-hli-1}
Let $\{\phi_n^{\tau}\}_{_{n,\tau}}=\Phi_0^{(\infty)}$ be as above. Then $\phi_n^{\tau}$ maps $\bigwedge^r H^1(k_n(\tau)_p,T)$ isomorphically onto $\al_n^{\tau} (=\textup{ker}(\Psi_n^{\tau}))$,  for all $\tau$ and $n$. In particular $\Phi_0^{(\infty)}$ satisfies $\hli$. 
\end{prop}

\begin{proof}
The proof is identical to the proof of (the easy half of) Lemma 3.1 of~\cite{kbb}, which also follows the proof of~\cite{mr02} Lemma B.1  almost line by line.
\end{proof}

\begin{rem}
\label{rem:lifts}
Note that $\Phi_0^{(\infty)}$ lifts, by construction, the element $\Phi_0$ of~\cite{kbb} \S2.3 in the sense we described above. 
\end{rem}
\begin{rem}\label{interchange}
Since the maps $$H^1(k_p,T\otimes\LL)=H^1_{\textup{Iw}}(k_p,T):=\varprojlim_n H^1((k_n)_{p},T) \lra H^1((k_n)_{p},T)$$ are all surjective and $H^1(k_p,T\otimes\LL)$ (\emph{resp.} $H^1((k_n)_{p},T)$) is a free $\LL$ (\emph{resp.} $\ZZ_p[\Gamma_n]$) module of rank $r$, it follows that there is a canonical isomorphism $$\varprojlim_n\bigwedge^r H^1((k_n)_p,T) \cong \bigwedge^r \varprojlim_n H^1((k_n)_p,T)= \bigwedge^r H^1(k_p,T\otimes\LL).$$
This, together with  Proposition~\ref{homs-hli-1}, shows that $\phi_\infty=\{\phi_n\}_{_n}$ maps $\bigwedge^r H^1(k_p,T\otimes\LL)$ isomorphically onto $\LLL=\varprojlim_n \al_n$.
\end{rem}

\subsection{Kolyvagin systems for the $\LLL$-modified Selmer triple-II}
\label{kolsys-2}
We now ready to construct  a $\LL$-adic Kolyvagin system\footnote{See Definitions 3.1.3 and 3.1.6 of~\cite{mr02} for a precise definition of a Kolyvagin system.} $\{\kappa_{\eta}^{\Phi}\}$ for the $\LLL$-modified Selmer structure $\FF_{\LLL}$ on $T\otimes\LL$, using the Euler system $\{\varepsilon_{k_n(\tau),\Phi}^{\chi}\}_{_{n,\tau}}$ of  $\Phi$-Stark elements, for each choice of a compatible homomorphisms $\Phi \in \hhh$ that satisfies $\hli$. We will use these classes in the next section to prove our main results.

Let $\PP$ denote the set primes of $k$ whose elements do not divide $pf_{\chi}$. For each positive integer $m$, let $$\PP_{m+n}=\{\qq \in \PP: \qq \hbox{ splits completely in } L(\mu_{p^{m+n+1}})/k\}$$ be a subset of $\PP$. Note that $\PP_j$ is exactly the set of primes being determined by Definition 3.1.6 of~\cite{mr02} or \S 4, Definition 1.1 of~\cite{r00} when $T=\ZZ_p(1)\otimes\chi^{-1}$. Let $\NN$ (\emph{resp.} $\NN_j$) denote the square free products of primes $\qq$ in $\PP$ (\emph{resp.} $\PP_j$). We also include $1$ in $\NN$ (and $\NN_j$). For notational simplicity, we also write $\mathbb{T}:=T\otimes\LL$, and for a fixed topological generator $\gamma$ of $\Gamma=\Gal(k_\infty/k)$ we set $\gamma_n=\gamma^{p^n}$.

Theorem 5.3.3 of~\cite{mr02} gives a map $$\textup{\textbf{ES}}(T)\lra \overline{\textup{\textbf{KS}}}(\mathbb{T},\FF_{cl},\PP)$$ where $\textup{\textbf{ES}}(T)$ denote the collection of Euler systems for $T$, and $$\overline{\textup{\textbf{KS}}}(\mathbb{T},\FF_{cl},\PP):=\lim_{\stackrel{\longleftarrow}{m,n}}(\lim_{\stackrel{\lra}{j}} \textup{\textbf{KS}}(\mathbb{T}/(p^m,\gamma_{n}-1)\mathbb{T},\,\FF_{cl},\PP \cap \PP_{j}))$$ such that if an Euler system $\{c_{k_n(\tau)}\}_{\tau}=\textbf{c} \in \textup{\textbf{ES}}(T)$ maps to the Kolyvagin system $\kappa=\{\kappa_{\eta}(m,n)\}_{\eta \in \NN_{m+n}}$ under this map, then 
\begin{align*}
H^1(\QQ,\mathbb{T})=\varprojlim_{m,n}H^1(\QQ,\mathbb{T}&/(p^m,\gamma_n-1)\mathbb{T}) \ni \varprojlim_{m,n} \kappa_{1}(m,n):=\kappa_1=\\ &=\{\textbf{c}_{k_n}\} \in \varprojlim_n H^1(\QQ,\mathbb{T}/(\gamma_n-1)\mathbb{T}) =H^1(\QQ,\mathbb{T})
\end{align*}

Let $\{\kappa^{\Phi}_{\eta}(m,n)\}_{_{\eta \in \NN;m,n}}$ be the image of the Euler system $\{\varepsilon_{k_n(\tau),\Phi}^{\chi}\}_{_{n,\tau}}$ of Proposition~\ref{def: es}.

\begin{rem}
\label{shapiro}
Using Shapiro's lemma one easily checks that $$H^1(k(\tau),\mathbb{T}/(p^m,\gamma_{n}-1)\mathbb{T}) \cong H^1(k_n(\tau),T/p^mT), \hbox{ and}$$ $$H^1(k(\tau)_p,\mathbb{T}/(p^m,\gamma_{n}-1)\mathbb{T}) \cong H^1(k_n(\tau)_p,T/p^mT).$$ See~\cite{mr02} Lemma 5.3.1 and~\cite{r00} Appendix B.5 for the semi-local version. We will make use of these identifications below.
\end{rem}
\begin{thm}
\label{mainkolsys}
Assume that $\Phi \in \hhh$ satisfies $\hli$. Then $$\kappa^{\Phi}:=\{\kappa^{\Phi}_{\eta}(m,n)\}_{_{\eta \in \NN;m,n}} \in  \overline{\textup{\textbf{KS}}}(\mathbb{T},\FF_{\LLL},\PP).$$
\end{thm}

Here $$\overline{\textup{\textbf{KS}}}(\mathbb{T},\FF_{\LLL},\PP)= \lim_{\stackrel{\longleftarrow}{m,n}}(\lim_{\stackrel{\lra}{j}} \textup{\textbf{KS}}(\mathbb{T}/(p^m,\gamma_{n}-1)\mathbb{T},\,\FF_{\LLL},\PP \cap \PP_{j}))$$is the module of $\LL$-adic Kolyvagin systems for the $\LLL$-modified Selmer structure $\FF_{\LLL}$ on $\mathbb{T}$. An element $\kappa^{\Phi} \in \overline{\textup{\textbf{KS}}}(\mathbb{T},\FF_{\LLL},\PP)$ will be called a \emph{$\LL$-adic Kolyvagin system of $\Phi$-Stark elements}. 

\begin{rem}
\label{1-kolsys}Using the fact that both $\{p^m,\gamma_n-1\}_{_{m,n}}$ and $\{p^m,\xx^n\}_{_{m,n}}$ form a base of neighborhoods at $0$, it is easy to see that the module of Kolyvagin systems $\overline{\textup{\textbf{KS}}}(\mathbb{T},\FF_{\LLL},\PP)$ defined as above is naturally isomorphic to the module of Kolyvagin systems $\overline{\textup{\textbf{KS}}}(T\otimes\LL,\FF_{\LLL})$ defined as in~\cite{kbb} \S\S3.1-3.2. By abuse of this observation, a $\LL$-adic Kolyvagin system of $\Phi$-Stark elements $\kappa^{\Phi}$ will be regarded as an element of the module of Kolyvagin systems $\overline{\textup{\textbf{KS}}}(T\otimes\LL,\FF_{\LLL})$ of~\cite{kbb}. We further note that, one may also identify these modules by what Mazur and Rubin call the generalized module of Kolyvagin systems, \emph{c.f.} \cite{mr02} Definition 3.1.6.  
 \end{rem}
For the rest of this section the integers $m$ and $n$ will be fixed, and we will simply denote the element $\kappa^{\Phi}_{\eta}(m,n) \in H^1(k,\mathbb{T}/(p^m,\gamma_{n}-1)\mathbb{T})$ by $\kappa^{\Phi}_{\eta}$. Note that the above statement of Theorem~\ref{mainkolsys} says that for each $\eta \in \NN_{m+n}$, $\kappa^{\Phi}_{\eta} \in H^1_{\FF_{\LLL}(\eta)}(k,\mathbb{T}/(p^m,\gamma_{n}-1)\mathbb{T})$, where $\FF_{\LLL}(\eta)$ is defined as in Example 2.1.8 of~\cite{mr02}. However, Theorem 5.3.3 of~\cite{mr02} already says that $\kappa^{\Phi}_{\eta} \in  H^1_{\FF_{cl}(\eta)}(k,\mathbb{T}/(p^m,\gamma_{n}-1)\mathbb{T})$; therefore to prove Theorem~\ref{mainkolsys} it suffices to prove the following (since $\FF_{\LLL}$ and $\FFc$ determine the same local conditions at places not dividing $p$):

\begin{prop}
\label{loc-p}
Let $$\xymatrix{\textup{loc}_p: &H^1(k,\mathbb{T}/(p^m,\gamma_{n}-1)\mathbb{T}) \ar[r] \ar[d]^{\cong}& H^1(k_p,\mathbb{T}/(p^m,\gamma_{n}-1)\mathbb{T})\ar[d]^{\cong}\\
& H^1(k_n,T/p^mT)\ar[r] & H^1((k_n)_p,T/p^mT)
}$$ be the localization map into the semi-local cohomology at $p$ (vertical isomorphisms follow from Shapiro's lemma, c.f. Remark~\ref{shapiro}). Then $$\textup{loc}_p (\kappa^{\Phi}_{\eta}) \in \al_n/p^m\al_n \subset H^1((k_n)_p,T/p^mT)$$
\end{prop}
We prove Proposition~\ref{loc-p} below. We remark that $\al_n/p^m\al_n$ is the propagation $H^1_{\FF_{\LLL}(\eta)}((k_n)_p,T/p^mT)$ of the induced local condition $H^1_{\FF_{\LLL}(\eta)}((k_n)_p,T):=H^1_{\FF_{\LLL}}(k_p,\mathbb{T}/(p^m,\gamma_{n}-1)\mathbb{T})=\al_n$ at $p$ (which is also obtained by propagating the Selmer structure $\FF_{\LLL}$ on $T\otimes\LL=\mathbb{T}$ to $\mathbb{T}/(p^m,\gamma_{n}-1)\mathbb{T}$. Let $$\{\tilde{\kappa}^{\Phi}_{\eta}(m,n) \in H^1((k_n)_p,T/p^mT) \}_ {\eta \in \NN_{m+n}}$$ be the collection that Definition 4.4.10 of~\cite{r00} associates to the Euler system $\{\varepsilon_{k_n(\tau),\Phi}^{\chi}\}_{_{n,\tau}}$. Here we write $\tilde{\kappa}^{\Phi}_{\eta}(m,n)$ for the class denoted by $\kappa_{[k_n,\eta,m]}$ in~\cite{r00}. Since we have already fixed $m$ and $n$ until the end of this section, we will safely drop $m$ and $n$ from the notation and denote $\tilde{\kappa}^{\Phi}_{\eta}(m,n)$ by $\tilde{\kappa}^{\Phi}_{\eta}$ if there is no danger for confusion. Note that Equation (33) in Appendix A of~\cite{mr02} relates these class to $\kappa^{\Phi}_{\eta}$.
\begin{lemma}
\label{reduction-1}
Assume $\textup{loc}_p(\tilde{\kappa}^{\Phi}_{\eta}) \in \al_n^{\tau}/p^m \al_n^{\tau}$ then $\textup{loc}_p(\kappa^{\Phi}_{\eta}) \in \al_n^{\tau}/p^m \al_n^{\tau}$ as well.
\end{lemma}
\begin{proof}
Obvious using Equation (33) in Appendix A of~\cite{mr02}.
\end{proof}

Let $D_{\eta}$ denote the derivative operators as in Definition 4.4.1 of~\cite{r00}. Definition 4.4.10 (and Remark 4.4.3) defines $\tilde{\kappa}^{\Phi}_{\eta}$ as the inverse image of $D_{\eta}\varepsilon_{k_n(\eta),\Phi}^{\chi}$ (mod~$p^m$) under the restriction map\footnote{ Note that $(\mu_{p^{\infty}}\otimes \chi^{-1})^{G_{k_n(\eta)}}$ is trivial (where $G_{k_n(\eta)}$ stands for the absolute Galois group of the totally real field $k_n(\eta)$),   since, for example, complex conjugation cannot act by $\chi$ on $\mu_{p^{\infty}}$ since $\chi$ is even. This argument proves that this restriction map is an isomorphism, by Remark 4.4.3 of~\cite{r00}.}
\begin{align*} H^1(k_n,T/p^mT) &\lra H^1(k_n(\eta),T/p^mT)^{G_{\eta}}\end{align*} Therefore $\textup{loc}_p(\tilde{\kappa}^{\Phi}_{\eta})$ maps to $\textup{loc}_p(D_{\eta}\varepsilon_{k_n(\eta),\Phi}^{\chi})$ (mod~$p^m$) under the map (which is also an isomorphism by Remark 4.4.3, Proposition B.5.1 and Proposition B.4.2 of~\cite{r00}) $$H^1((k_n)_p,T/p^mT) \lra H^1(k_n(\eta)_p,T/p^mT)^{G_{\eta}}$$ 
Under this isomorphism $\al_n/p^m\al_n$ is mapped isomorphically onto the $\ZZ/p^m\ZZ[G_\eta \times\Gamma_n]$-line  $[\al_n^{\eta}/p^m\al_n^{\eta}]^{G_{\eta}}$, by definition of $\al_n^{\eta}$ and the fact that it is a free $\ZZ_p[G_{\eta}\times\Gamma_n]$-module . Below diagram summarizes the discussion in this paragraph:
$$\xymatrix{
H^1(k_p,T/p^mT) \ar[r]^(.45){\sim}& H^1(k(\eta)_p,T/p^mT)^{G_{\eta}}\\
\al_n^{\tau}/p^m\al_n^{\tau}\ar[r]^(.4){\sim} \ar@{^{(}->}[u] &[\al_n^{\eta}/p^m\al_n^{\eta}]^{G_{\eta}}\ar@{^{(}->}[u] 
}$$

\begin{prop}
\label{lined}
Suppose $\Phi$ satisfies $\hli$. Then $\textup{loc}_p(\tilde{\kappa}^{\Phi}_{\eta}) \in \al_n^{\eta}/p^m \al_n^{\eta}$.
\end{prop}

\begin{proof}
Since $\textup{loc}_p$ is Galois equivariant $\textup{loc}_p(D_{\eta}\varepsilon_{k_n(\eta),\Phi}^{\chi})=D_{\eta}\textup{loc}_p(\varepsilon_{k_n(\eta),\Phi}^{\chi})$. Further $\textup{loc}_p(\varepsilon_{k_n(\eta),\Phi}^{\chi}) \in \al_n^{\eta}$ since $\Phi$ satisfies $\hli$. On the other hand, by Lemma 4.4.2 of~\cite{r00} $D_{\eta}\varepsilon_{k_n(\eta),\Phi}$ (mod~$p^m$) is fixed by $G_{\eta}$, which in return implies $$\textup{loc}_p(\varepsilon_{k_n(\eta),\Phi}^{\chi}) \,(\textup{mod\,} p^m) \in [\al_n^{\eta}/p^m\al_n^{\eta}]^{G_{\eta}}$$ This proves that $\textup{loc}_p(\tilde{\kappa}^{\Phi}_{\eta})$ maps into $\al_n/p^m \al_n$ by above discussion.
\end{proof}
\begin{proof}[Proof of Proposition~\ref{loc-p}]
Immediately follows from Lemma~\ref{reduction-1} and Proposition~\ref{lined}.
\end{proof}
By the discussion following the statement of Theorem~\ref{mainkolsys}, this also completes the proof of Theorem~\ref{mainkolsys}.

\newpage
\section{Applications to main conjectures}
\label{applications}
Assume throughout this section that the set of places $S$ of $k$ (which was introduced in \S\ref{sec:stark-euler}) contains no non-archimedean prime of $k$ which splits completely in $L/k$. We also assume that Leopoldt's conjecture holds true. Let $\Phi_0^{(\infty)} \in \hhh$ be as in \S\ref{sec:homs} and let $\kappa^{\Phi_0^{(\infty)}}=\{\kappa_\eta^{\Phi_0^{(\infty})}\} \in \overline{\KS}(\mathbb{T},\FF_{\LLL},\PP)$ be the $\LL$-adic Kolyvagin system of $\Phi_0^{(\infty)}$-Stark elements (which were introduced in \S\ref{kolsys-2}).

In~\cite{kbbstark}, under our hypotheses,  we explicitly determine a generator for the $\ZZ_p$-module of Kolyvagin systems $\overline{\KS}(T,\FF_{\mathcal{L}})$ for the Selmer structure $\FF_\al$ on $T$ using the Euler systems of Stark elements. (This module is known to be free of rank one as a $\ZZ_p$-module by~\cite{mr02} Theorem 5.2.10, and to be generated by a \emph{primitive} Kolyvagin system, in the sense of~\cite{mr02} Definition 4.5.5)  We denote this generator $\kappa^{\Phi_0}=\{\kappa_\eta^{\Phi_0}\}$ following~\cite{kbbstark}, where  $\Phi_0$ is an element of $\varprojlim_\tau \bigwedge^{r-1}\hbox{Hom}_{\ZZ_p[G_\tau]}(H^1(k(\tau)_p,T), \ZZ_p[G_\tau])$ which was already constructed in~\cite{kbbstark} in an identical manner to the construction of $\Phi_0^{(\infty)}$.

It is now clear that $\kappa^{\Phi_0^{(\infty)}}$ maps to the element $\kappa^{\Phi_0}$ under the map (which is surjective, by Theorem~\ref{main-stark} and Remark~\ref{1-kolsys})
$$\overline{\textup{\textbf{KS}}}(\mathbb{T},\FF_{\LLL},\PP)  \lra\overline{\KS}(T,\FF_{\mathcal{L}}). $$
On the other hand, since the Kolyvagin system  $\kappa^{\Phi_0}$ is primitive, it follows that 
\begin{cor}
\label{primitive}
The $\LL$-adic Kolyvagin system of $\Phi_0^{(\infty)}$-Stark elements $\kappa^{\Phi_0^{(\infty)}}$ is $\LL$-primitive, in the sense of~\cite{mr02} 5.3.9.
\end{cor}

Let $\textup{char}(\mathbb{A})$ denote the characteristic ideal of a torsion $\LL$-module $\mathbb{A}$. Recall that $A^{\vee}$ denotes the Pontryagin dual of an abelian group $A$. Recall also the definition of the dual Selmer module. The main application of the ($\LL$-primitive) Kolyvagin system of $\Phi_0^{(\infty)}$-Stark elements is the following:
\begin{thm}
\label{main-application}
$$\textup{char}\left(H^1_{\FF_{\LLL}}(k,\mathbb{T})/\LL\cdot\kappa^{\Phi_0^{(\infty)}}_1\right)=\textup{char}\left( (H^1_{\FF_{\LLL}^*}(k,\mathbb{T}^*)^{\vee}\right)$$
\end{thm}

\begin{proof}
This is~\cite{mr02} Theorem 5.3.10(iii) applied to our setting.
\end{proof}
  On the other  hand, Proposition~\ref{exact-seq-stark} shows immediately that
\begin{prop}\label{exact-seq-stark-2}
The sequence \[
\begin{array}{rl}
  0&\lra H^1_{\FF_{\LLL}}(k,\mathbb{T})/\LL\cdot\kappa^{\Phi_0^{(\infty)}}_1\lra \LLL/\LL\cdot\kappa^{\Phi_0^{(\infty)}}_1 \lra   \\\\
 &H^1_{\FF_{\textup{str}}^*}(k,\mathbb{T}^*)^{\vee} \lra (H^1_{\FF_{\LLL}^*}(k,(\mathbb{T})^*)^{\vee}\lra 0   \\
\end{array}
\]  is exact. 
\end{prop}

\begin{cor}
\label{first reduction}
$$\textup{char}\left(\LLL/\LL\cdot\kappa^{\Phi_0^{(\infty)}}_1\right)=\textup{char}\left((H^1_{\FF_{\textup{str}}^*}(k,\mathbb{T}^*)^{\vee}\right).$$
\end{cor}
\begin{proof}
Follows from Theorem~\ref{main-application} and Proposition~\ref{exact-seq-stark-2}.
\end{proof}
 
 Let $\mathbf{c}_{k_{\infty}}^{\textup{stark}} \in \bigwedge^rH^1(k_p,\mathbb{T})$ be the image of the element $$\{\varepsilon_{k_n}^{\chi}\}_{_n}\in \varprojlim_n \bigwedge^r H^1(k_n,T)\stackrel{\textup{loc}_p}{\lra}  \varprojlim_n \bigwedge^r H^1((k_n)_p,T) \stackrel{(\ref{interchange})}{\cong}H^1(k_p,T\otimes\LL)$$ under the localization at $p$ map $\textup{loc}_p$ composed with the isomorphisms of Remark~\ref{interchange} (we recall that $\mathbb{T}:=T\otimes\LL$).
 Further, again by Remark~\ref{interchange}, the element $\Phi_0^{(\infty)}$ maps $\bigwedge^rH^1(k_p,\mathbb{T})$ isomorphically onto $\LLL$ (and maps the element $\mathbf{c}_{k_{\infty}}^{\textup{stark}} \in \bigwedge^rH^1(k_p,\mathbb{T})$ to $\kappa^{\Phi_0^{(\infty)}}_1 \in \LLL$), and thus it follows that $$\bigwedge^rH^1(k_p,\mathbb{T})/\LL\cdot\mathbf{c}_{k_{\infty}}^{\textup{stark}} \cong \LLL/\LL\cdot\kappa^{\Phi_0^{(\infty)}}_1$$ Therefore we proved: 
\begin{cor}
\label{second reduction}
$$\textup{char}\left(\bigwedge^rH^1(k_p,\mathbb{T})/\LL\cdot\mathbf{c}_{k_{\infty}}^{\textup{stark}}\right)=\textup{char}\left(H^1_{\FF_{\textup{str}}^*}(k,\mathbb{T}^*)^{\vee}\right).$$
\end{cor}

Let $\mathcal{L}_{k}^{\chi}$ denote the Deligne-Ribet $p$-adic $L$-function attached to the character $\chi$ (see~\cite{deligne-ribet} for the construction of this $p$-adic $L$-function). In the spirit of~ \cite{pr-kubota}, \cite{pr2} and the explicit reciprocity laws conjectured by Perrin-Riou~\cite{pr} (and proved by  Colmez~\cite{colmez-reciprocity}) we propose the following:
\begin{conjecture}
\label{conj-kbb}
$\mathcal{L}_{k}^{\chi}$ generates $\textup{char}\left(\bigwedge^rH^1(k_p,\mathbb{T})/\LL\cdot\mathbf{c}_{k_{\infty}}^{\textup{stark}}\right)$.
\end{conjecture}

Using the explicit description of the Galois cohomology groups in question (see for instance~\cite{r00} 1.6.B) one may identify  $H^1_{\FF_{\textup{str}}^*}(k,\mathbb{T}^*)^{\vee}$ by $\Gal(M_\infty/L_\infty)^{\chi}$, where $M_\infty$ is the maximal abelian $p$-extension of $L_\infty$ unramified outside primes above $p$. This is the Iwasawa module in this setting. Hence, by Corollary~\ref{second reduction}, Conjecture~\ref{conj-kbb} is equivalent to Iwasawa's main conjecture over the totally real field $k$ for characters $\chi$ with the properties ascribed at the beginning of  this paper. 

When $k=\QQ$ (so $r=1$) Conjecture~\ref{conj-kbb} is a Theorem of Iwasawa~\cite{iwasawa-1}. In fact in this case Stark elements are obtained from cyclotomic units. The key observation is that cyclotomic units both demonstrate the complex Stark conjecture and the $p$-adic Stark conjecture simultaneously. We do not think one could prove Conjecture~\ref{conj-kbb} directly using \emph{only} Rubin's Stark elements (which are solutions to Stark conjectures for complex $L$-functions), in fact one should rather use the solutions to an appropriate $p$-adic Stark conjecture. 

Having said that, we note that this conjecture is in fact true by the work of Wiles~\cite{wiles-mainconj}. This remark already points at a relation between the solutions of complex Stark conjectures and $p$-adic Stark conjectures. This relation should be understood as an analogy to the fact that cyclotomic units give solutions to complex and $p$-adic Stark conjectures (when $k=\QQ$) simultaneously. 

In a sequel to this paper, we hope to discuss the relation between solutions to complex and $p$-adic Stark conjectures via the rigidity offered by Theorem~\ref{main-stark}, and prove Conjecture~\ref{conj-kbb} (without assuming Wiles' work) thus deduce the main conjectures. 

 \newpage
\appendix
\section{Local Conditions at $p$ over an Iwasawa algebra via $(\phi,\Gamma)$-modules}
\label{fontaine}

Throughout this section let $K$ denote a finite extension of $\QQ_p$ and set $\tilde{K}_n:=K(\mu_{p^n})$, $\tk_{\infty}:=\bigcup_{n}\tk_n$.  Set also $\tilde{H}_K:=\textup{Gal}(\overline{K}/\tk_{\infty})$ and $\tGamma_K:=G_K/\tH_K=\textup{Gal}(\tk_{\infty}/K)$. Let $\tgamma$ be a topological generator of the pro-cyclic group $\tGamma_K$, and let $\tlambda_K:=\ZZ_p[[\tGamma_K]]$. Let $\tgamma_n$ be a fixed topological generator of $\textup{Gal}(\tk_{\infty}/\tk_n):=\tGamma_n$ for $n \in \ZZ^+$, chosen in a way that $\tgamma_n=\tgamma_1^{p^{\alpha_n}}$, where $\alpha_n \in \ZZ^+$ (so that $[\tk_n:\tk]=p^{\alpha_n}$).

 Let $K_n$ be the maximal pro-$p$ extension of $K$ inside $\tk_n$, and let $\bigcup_n K_n=:K_{\infty} \subset \tk_{\infty}$ be the cyclotomic $\ZZ_p$-extension of $K$. We set $\Gamma_K:=\textup{Gal}(K_{\infty}/K)$ and $\Lambda_K:=\ZZ_p[[\Gamma_K]]$. Note then that $$\tGamma_K=W\times \Gamma_K,\hbox{ and that } \tlambda_K=\ZZ_p[W]\otimes_{\ZZ_p}\Lambda_K$$ where $W$ is a finite group which has prime-to-$p$ order (In fact $W$ can be identified by $\textup{Gal}(K(\mu_p)/K)$). Let $\gamma$ denote the restriction of $\tgamma$ to $K_{\infty}$; so that $\gamma$ is a topological generator of $\Gamma_K$. Let also $\gamma_n$ denote the image of $\tgamma_n$ under the natural isomorphism $$\textup{Gal}(\tk_{\infty}/\tk_n) \cong \textup{Gal}(K_{\infty}/K)$$ and set $H_K:=\textup{Gal}(\overline{K}/K_{\infty})$ (so that $H_K/\tH_K\cong W$).
 
 In~\cite{fontaine-phi-gamma}, \footnote{However, one should be cautious as Fontaine uses $K_n$ for our $\tk_n$, etc. For example, his $\Gamma_K$ is our $\tGamma_K$.}Fontaine introduces the notion of a $(\phi,\Gamma)$-module over certain period rings that he denotes by $\fontaineO$ (which really is the ring of integers of the field $\widehat{\varepsilon^{\textup{nr}}}=\textup{Frac}(\fontaineO))$; we also set $\mathcal{O}_{\varepsilon(K)}:=(\mathcal{O}_{\widehat{\varepsilon^{\textup{nr}}}})^{H_K}$. We will not include a  detailed discussion of these objects here, we refer the reader to A.3.1-3.2 of~\cite{fontaine-phi-gamma} for their definitions and basic properties. Briefly, a $(\phi,\Gamma)$-module over $\fontaineOK$ is a finitely generated $\fontaineOK$-module with semi-linear continuous and commuting actions of $\phi$ and $\Gamma:=\tGamma_K$. A $(\phi,\Gamma)$-module $D$ over $\fontaineOK$ is \emph{\'etale} if $\phi(D)$ generates $D$ as an $\fontaineOK$-module.
 
 Using his theory, Fontaine established an equivalence between the category of $\ZZ_p$-representations of $G_K$ and the category of \'etale $(\phi,\Gamma)$-modules over $\fontaineOK$. The equivalence is given by 
 $$\begin{array}{rcl}
 T&\longmapsto &D(T):=(\fontaineO \otimes_{\ZZ_p}T)^{H_K}\\
 (\fontaineO \otimes_{\fontaineOK}D)^{\phi=1}=:T(D)&\longmapsfrom&D
 \end{array}
 $$ See A.1.2.4-1.2.6 in~\cite{fontaine-phi-gamma}.
 
 Herr~\cite{herr-1} later showed that one may use $(\phi,\Gamma)$-modules to compute the Galois cohomology groups $H^*(K,T)$. One of the benefits of his approach is that the complex he constructs with cohomology $H^*(K,T)$ is quite explicit, which allows us to compute certain local Galois cohomology groups of $p$-adic fields. In~\cite{herr-2}, he gives a proof of local Tate duality, where again the pairing is fairly explicitly constructed (see \S5 of~\cite{herr-2}) using the residues of the differential forms on $\fontaineOK$. The rest of this section is a survey of his result, and of its applications to Iwasawa theory following~\cite{benois} and~\cite{colmez}.

In the theory of $(\phi,\Gamma)$-modules, there is an important operator\footnote{ This definition makes sense because: \begin{itemize}
\item  $\textup{Tr}_{\widehat{\varepsilon^{\textup{nr}}}/\phi(\widehat{\varepsilon^{\textup{nr}}})} \subset p\fontaineO$,
\item $\phi$ is injective.
\end{itemize}
} $$\psi: \fontaineO \lra \fontaineO$$ $$\psi(x):= \frac{1}{p}\phi^{-1}(\textup{Tr}_{\widehat{\varepsilon^{\textup{nr}}}/\phi(\widehat{\varepsilon^{\textup{nr}}})}(x))$$ which will be crucial for what follows. It is a left inverse of $\phi$, and its action on $\fontaineO$ commutes with the action of $G_K$. It induces an operator $$\psi: D(T) \lra D(T)$$ for any $G_K$-representation $T$. 

Let $C_{\psi,\gamma}$ be the complex  $$\xymatrix{
C_{\psi,\gamma}: \,0 \ar[r] &D(T) \ar[rr]^(.4){(\psi-1,\gamma-1)} &&D(T) \oplus D(T) \ar[rr]^(.57){(\gamma-1) \ominus (\psi-1)}&&D(T)\ar[r] &0
}$$ Main result of~\cite{herr-1} is:

\begin{Athm}
\label{herr-main} 
The complex $C_{\psi,\gamma}$ computes the $G_K$-cohomology of $T$:
\begin{enumerate}
\item[\textbf{(i)}] $H^0(K,T)\cong D(T)^{\psi=1, \tgamma=1}$,
\item[\textbf{(ii)}]$H^2(K,T) \cong D(T)/(\psi-1,\tgamma-1)$,
\item[\textbf{(iii)}] There is an exact sequence 
$$0 \lra \frac{D(T)^{\psi=1}}{\tgamma-1} \lra H^1(K,T)\lra \left(\frac{D(T)}{\psi-1}\right)^{\tgamma=1} \lra 0$$
\end{enumerate}
Further, all the isomorphisms and maps that appear above are functorial in $T$ and $K$.
\end{Athm}

\begin{Adefine}
\label{iwasawa-p}
Let $$H^1_{\tilde{\textup{Iw}}}(K,T):=\varprojlim_{n}(H^1(\tk_n,T)) \hbox{ and }$$ $$H^1_{\textup{Iw}}(K_n,T):=\varprojlim_{n}(H^1(K_n,T))$$ where the inverse limits are with respect to the corestriction maps. 
\end{Adefine}

\begin{Arem}
\label{K-tilde-vs-K}
Since the order of $W$ is prime to $p$, it follows that $$H^1_{\textup{Iw}}(K_n,T) \stackrel{\sim}{\lra}H^1_{\tilde{\textup{Iw}}}(K,T)^W$$ by Hochschild-Serre spectral sequence.
\end{Arem}
We now determine the structure of $H^1_{\textup{Iw}}(K_n,T)$ using Theorem~A.\ref{herr-main}.
\begin{Aprop}
\label{nth-layer-p}
Define $\tau_n:=1+\tgamma_{n-1}+ \dots + \tgamma_{n-1}^{p-1} \in \ZZ_p[[\tGamma_K]]$. Then we have a commutative diagram with exact rows:
$$\xymatrix{C_{\psi,\tgamma_n}(\tk_n, T): 0 \ar @{.>}[d]^{\tau_n^{*}} \ar[r] &D(T) \ar[r]\ar[d]^{\tau_n}& D(T) \oplus D(T) \ar[r] \ar[d]^{\tau_n \oplus \textup{id}}& D(T) \ar[r]\ar[d]^{\textup{id}}& 0\\
C_{\psi,\tgamma_{n-1}}(\tk_{n-1}, T): 0 \ar[r] &D(T) \ar[r]& D(T) \oplus D(T) \ar[r]& D(T) \ar[r]& 0
}$$
and the map induced from the morphism $\tau_n^*$ on the cohomology of $C_{\psi,\tgamma_n}(\tk_n, T)$ is the corestriction map under Herr's identification (Theorem~A.\ref{herr-main}) $$H^*(C_{\psi,\tgamma_n}(\tk_n, T)) \cong H^*(\tk_n,T).$$
\end{Aprop}

\begin{proof}
$\tau_n^*$ is a cohomological functor and induces $\textup{Tr}_{\tk_n/\tk_{n-1}}$ on $H^0$, so it induces corestrictions on $H^i$.
\end{proof}

Using Proposition~A.\ref{nth-layer-p} one may compute $H^*_{\tilde{\textup{Iw}}}(K,T)$:
\begin{Athm}
\label{cohomology-iwasawa-p-tilde}
\begin{enumerate}
\item[\textbf{(i)}] $H^i_{\tilde{\textup{Iw}}}(K,T)=0$ if $i \neq 1,2$,
\item[\textbf{(ii)}] $H^1_{\tilde{\textup{Iw}}}(K,T)\stackrel{\sim}{\lra} D(T)^{\psi=1}$,
\item[\textbf{(iii)}] $H^2_{\tilde{\textup{Iw}}}(K,T) \stackrel{\sim}{\lra} D(T)/(\psi-1)$.
\end{enumerate}

\end{Athm}
See \S{II.3} of~\cite{colmez} for a proof of this theorem.
\begin{Arem}
\label{rem:p-adic L-function}
The isomorphism $$\exp^*: H^1_{\tilde{\textup{Iw}}}(K,T)\stackrel{\sim}{\lra} D(T)^{\psi=1}$$ of Theorem~A.\ref{cohomology-iwasawa-p-tilde} is the map that will give rise to the conjectural $p$-adic $L$-function.
\end{Arem}

Let $\mathcal{C}(T):=(\phi-1)D(T)^{\psi=1}$. Since $\psi$ is a left inverse of $\phi$, it follows that $$\ker\{D(T)^{\psi=1} \lra \mathcal{C}(T)\}=D(T)^{\phi=1}$$ hence we have an exact sequence:
\be
\label{phi-psi-sequence}
0 \lra D(T)^{\phi=1} \lra D(T)^{\psi=1} \stackrel{\phi-1}{\lra}\mathcal{C}(T) \lra 0
\ee 

Using techniques from the theory of $(\phi,\Gamma)$-modules one can determine the structure of $\mathcal{C}(T)$:
\begin{Aprop}
\label{free-C}
$\mathcal{C}(T)$ is free of rank $[K:\QQ_p] \cdot\textup{rank}_{\ZZ_p}T$ as a $\ZZ_p[[\tGamma_K]]$-module. 
\end{Aprop}

Also one may check that $D(T)^{\phi=1}\cong T^{\tH_K}$, which is finitely generated over $\ZZ_p$, hence is a torsion $\ZZ_p[[\tGamma_K]]$-module. Thus, it follows from Proposition~A.\ref{free-C} and Theorem~A.\ref{cohomology-iwasawa-p-tilde} that $D(T)^{\phi=1}=H^1_{\tilde{\textup{Iw}}}(K,T)_{\textup{tors}}$, the torsion submodule of $H^1_{\tilde{\textup{Iw}}}(K,T)$.

If we now take $W$-invariance of the exact sequence~\ref{phi-psi-sequence}, use the fact that taking $W$-invariance is an exact functor; use Theorem~A.\ref{cohomology-iwasawa-p}, Remark~A.\ref{K-tilde-vs-K} and Proposition~A.\ref{free-C} we obtain:

\begin{Athm}\textup{(see~\cite{colmez})}
\label{cohomology-iwasawa-p}
\begin{enumerate}
\item[\textbf{(i)}] $H^1_{\textup{Iw}}(K,T)_{\textup{tors}} \cong T^{H_K}$, 
\item[\textbf{(ii)}] The  $\ZZ_p[[\Gamma_K]]$-module $H^1_{\textup{Iw}}(K,T)/H^1_{\textup{Iw}}(K,T)_{\textup{tors}}$ is free of rank $[K:\QQ_p]\cdot \textup{rank}_{\ZZ_p}T$.
\end{enumerate}
\end{Athm}

\bibliographystyle{alpha}
\bibliography{biblio-tez}


\end{document}